\documentclass[12pt]{article}
\usepackage[utf8]{inputenc}

\usepackage{amsmath}
\usepackage{amsfonts}
\usepackage{amssymb}
\usepackage{anysize}
\usepackage{cite}
\usepackage{comment}
\usepackage[shortlabels]{enumitem}
\usepackage{mathrsfs}
\usepackage{amsthm}
\usepackage{tikz}
    \usetikzlibrary{matrix}
    \usetikzlibrary{decorations.markings}
\usepackage{hyperref}
\usepackage[T1]{fontenc}
\usepackage{oldgerm}
\usepackage{scalerel,stackengine}
\usepackage{stmaryrd}
\usepackage{xcolor}
\usepackage{xfrac}
\usepackage{authblk}

\stackMath
\newcommand\reallywidehat[1]{%
\savestack{\tmpbox}{\stretchto{%
  \scaleto{%
    \scalerel*[\widthof{\ensuremath{#1}}]{\kern.1pt\mathchar"0362\kern.1pt}%
    {\rule{0ex}{\textheight}}%WIDTH-LIMITED CIRCUMFLEX
  }{\textheight}% 
}{2.1ex}}%
\stackon[-6.9pt]{#1}{\tmpbox}%
}
\newcommand\varfrak[1]{\mathord{\text{\textgoth{#1}}}}

\newtheorem{theorem}{Theorem}[section]
\newtheorem{definition}{Definition}[section]
\newtheorem{example}{Example}[section]
\newtheorem{proposition}{Proposition}[section]
\newtheorem{lemma}{Lemma}[section]
\newtheorem{corollary}{Corollary}[section]
\newtheorem{remark}{Remark}[section]

\begin{document}

\title{Principal bundles and connections modelled by Lie group bundles}

\author[1]{Marco Castrillón López\thanks{mcastri@mat.ucm.es}}

\author[2]{Álvaro Rodríguez Abella\thanks{alvrod06@ucm.es}}

\affil[1]{Facultad de Ciencias Matemáticas, Universidad Complutense de Madrid Plaza de las Ciencias, 3,Madrid, 28040,Madrid, Spain}

\affil[2]{Insituto de Ciencias Matemáticas (CSIC-UAM-UC3M-UCM), Calle Nicolás Cabrera, 13-15, Madrid, 28049, Madrid, Spain}

%\affil[$\spadesuit$]{\normalsize These authors contributed equally to this work.}

\date{}

\maketitle

\noindent \emph{MSC2020:} 53C05, 53C15, 22E99

\noindent \emph{Key words:} Ehresmann connection, fiberwise action, Lie group bundle, parallel transport, generalized principal bundle, Utiyama Theorem

\begin{abstract}
In this work, generalized principal bundles modelled by Lie group bundle actions are investigated. In particular, the definition of equivariant connections in these bundles, associated to Lie group bundle connections, is provided, together with the analysis of their existence and their main properties. The final part gives some examples. In particular, since this research was initially originated by some problems on geometric reduction of gauge field theories, we revisit the classical Utiyama Theorem from the perspective investigated in the article.
\end{abstract}

%\tableofcontents

\section{Introduction}

A fiberwise action of a Lie group fiber bundle on a smooth bundle over the same base leads to the notion of generalized principal bundle. They share some similar properties to (standard) principal bundles \cite[Ch.II]{kobayashinomizu1963}, but they also show important differences. For instance, on one hand it is possible to build (generalized) associated bundles by the action of Lie group bundles on other bundles. But on the other hand, Ehresmann connections \cite{saunders,michor1993} that are equivariant with respect to the action, require a precise and correct approach providing the corresponding notion of generalized principal connections. Unlike usual principal connections, they are always associated to a certain connection on the Lie group bundle, which gives an additional term in the equivariance formula. 

Lie group bundles, together with their natural infinitesimal companions, the Lie algebra bundles, are classical objects in the literature (see, for example, \cite{Douady1966}). These type of bundles naturally arise in various geometric contexts and in different applications. The reader should be aware that these bundles may involve questions concerning non-isomorphic algebraic structure between different fibers (the seminal work of Douady and Lazard \cite{Douady1966} already discussed situation; see \cite{Kira2019} for a recent approach). Although this scenario is very interesting and beautiful, in this work we confine ourselves to the case where the fibers are algebraically isomorphic (that is, the bundle is locally trivial from an algebraic point of view), a decision mainly motivated by the bundles that one encounters in the applications. With respect to them, they appear in a natural way when performing reduction by local symmetries in Lagrangian field theory \cite{abella2021,forger2012} and, since principal connections may be regarded from this perspective, they are also useful for classical reduction --that is, reduction by global symmetries-- in mechanics \cite{marsdenjurgen1993,cendramarsdenratiu2001,cendramarsdenratiu2001b} and field theories \cite{castrillonperezratiu2001,castrillonratiu2003,ratiu2011}. Lie group bundles are the unavoidable starting point in the geometric foundations gauge theories. In particular, they provide the basic language for a theory of geometric reduction of field theories when the group of symmetries are sections of Lie group bundles. This theory (that is still in progress) will collect the main instances of gauge theories and will require a wise use of the concepts exposed in this article. Actually, our initial interest in generalized principal bundles started in that framework, from where we have taken much inspiration. We would like to mention that fibered actions can be extended to the Lie groupoid setting (for example see \cite{Mackenzie2005,Hoyo2018,forger2015}), but we do not address this matter here.

Despite the interest and ubiquitous presence of these objects, it is remarkable to check the existence of some important gaps in the literature about the main properties, definitions, and the key geometric objects involved. In this work we aim at solving this situation with the study fibered actions, as well as the smooth bundles arising from the quotient by these actions: generalized principal bundles. Furthermore, we define equivariant connections on these bundles (i.e., generalized principal connections) and their curvature. Before doing that, we need to define Lie group bundle connections, which are connections on a Lie group bundle that respect the multiplicative structure. Actually, the generalized principal connections will be associated to Lie group bundle connections, a situation that does not have a counterpart in the notion of (standard) principal bundle connections.

The paper is organized as follows. In Section 2 we investigate fibered actions and quotients by them. After that, the definition of infinitesimal generators is recalled and generalized associated bundles are defined. In Section 3 we introduce Lie group bundle connections, as well as the induced linear connection on the Lie algebra bundle. Then generalized principal connections are defined and characterized using parallel transport. Besides, we prove a theorem of existence and study their curvature. In Section 4 we present several examples to illustrate the ideas of this work. In particular, we show that usual principal bundles and connections are particular cases of the generalized objects. We have a similar situation with connections in affine bundles. Finally, the action of the gauge group on connections is modelled with Lie group bundle actions. In this case, the generalized principal bundle provides a new approach to the well-known Utiyama Theorem. 

In the following, every manifold or map is smooth, meaning $C^\infty$. In addition, every fiber bundle $\pi_{Y,X}\colon Y\rightarrow X$ is assumed to be locally trivial and is denoted by $\pi_{Y,X}$. Given $x\in X$, $Y_x=\pi_{Y,X}^{-1}(\{x\})$ denotes the fiber over $x$. The space of (smooth) global sections of $\pi_{Y,X}$ is denoted by $\Gamma(\pi_{Y,X})$. In particular, vector fields on a manifold $X$ are denoted by $\mathfrak X(X)=\Gamma(\pi_{TX,X})$, where $TX$ is the tangent bundle of $X$. Likewise, the space of local sections on an open set $\mathcal U\subset X$ is denoted by $\Gamma(\mathcal U,\pi_{Y,X})$. The derivative, or tangent map, of a map $f\in C^\infty(X,X')$ between the manifolds $X$ and $X'$ is denoted by $(df)_x\colon  T_xX\rightarrow T_{x'}X'$, $x'=f(x)$. When working in local coordinates, we will assume the Einstein summation convention for repeated indices. A compact interval will be denoted by $I=[a,b]$.

%%%%%%%%%%%%%%%%%%%%
\section{Generalized principal bundles}

%%%%%%%%%%
\subsection{Actions of Lie group bundles}

A \emph{Lie group fiber bundle} with typical fiber a Lie group $G$ is a fiber bundle $\pi_{\mathcal G,X}:\mathcal{G} \to X$ such that for any point $x\in X$ the fiber $\mathcal G_x$ is equipped with a Lie group structure and there is a neighborhood $\mathcal U\subset X$ and a diffeomorphism $x\in\mathcal U\times G  \rightarrow \pi_{\mathcal G,X}^{-1}(\mathcal U)$ preserving the Lie group structure fiberwisely.

Note that the multiplication map $M\colon\mathcal G\times_X\mathcal G\to\mathcal G$ and the inversion map $\cdot^{-1}\colon\mathcal G\to\mathcal G$ are bundle morphisms covering the identity $\textrm{\normalfont id}_X\colon X\to X$, where $\times_X$ denotes the fibered product. Likewise, the map $1\colon X\rightarrow\mathcal G$ that assigns the identity element $1_x\in\mathcal G_x$ to each $x\in X$ is a global section (called the \emph{unit section}) of $\pi_{\mathcal G,X}$. Any Lie group bundle defines a Lie algebra bundle $\pi_{\varfrak g,X}\colon\varfrak g\to X$ as the vector bundle whose fiber $\varfrak g_x$ at each $x\in X$ is the Lie algebra of $\mathcal G_x$. That is, the Lie algebra bundle is the pull-back bundle $\varfrak g=1^*(V\mathcal G)$, where $V\mathcal G \subset T\mathcal G$ is the vertical bundle of $\pi_{\mathcal G,X}$, i.e. the kernel of $(\pi_{\mathcal G,X})_*$. 

\begin{remark}[Jets of Lie group fiber bundles]\label{example:jetliebundles}
Let $\pi_{\mathcal G,X}\colon\mathcal G\to X$ be a Lie group fiber bundle and $r\geq 0$ be an integer. Then the $r$-th jet bundle of $\pi_{\mathcal G,X}$, $J^r\mathcal G\to X$, is again a Lie group fiber bundle (see, for example, \cite[\S3, Th. 1]{forger2012}). The multiplication is inherited from the Lie group bundle structure of $\pi_{\mathcal G,X}$, that is,
\begin{equation*}
M\left(j_x^r\gamma_1,j_x^r\gamma_2\right)=j_x^r\left(M\circ(\gamma_1,\gamma_2)\right),\qquad j_x^r\gamma_1,j_x^r\gamma_2\in J^r\mathcal G.
\end{equation*}
\end{remark}

We consider subgroups of Lie group bundles in the following sense.

\begin{definition}
A \emph{Lie group subbundle} of a Lie group bundle $\pi_{\mathcal{G},X}:\mathcal{G}\to X$ is a Lie group bundle $\pi_{\mathcal H,X}\colon\mathcal H\to X$ such that $\mathcal H$ is a submanifold of $\mathcal{G}$ and $\mathcal H_x$ is a Lie subgroup of $\mathcal{G}_x$ for each $x\in X$. It is said to be \emph{closed} if $\mathcal H_x$ is a closed Lie subgroup of $\mathcal G_x$ for every $x\in X$. 
\end{definition}

Let $\pi_{Y,X}$ be a fiber bundle and $\pi_{\mathcal G,X}$ be a Lie group fiber bundle.

\begin{definition}
A \emph{(right) fibered action} of $\pi_{\mathcal G,X}$ on $\pi_{Y,X}$ is a bundle morphism \begin{equation*}
\Phi\colon Y\times_X\mathcal G\longrightarrow Y    
\end{equation*}
covering the identity $\mathrm{id}_X\colon X\to X$ such that $\Phi(y,hg)=\Phi(\Phi(y,h),g)$ and $\Phi(y,1_x)=y$, for all $(y,g),(y,h)\in Y\times_X\mathcal G$, $\pi_{\mathcal G,X}(y)=x$.
\end{definition}

For the sake of simplicity, we will denote $\Phi(y,g)=y\cdot g$ and we will say that $\pi_{\mathcal G,X}$ acts fiberwisely on the right on $\pi_{Y,X}$. Note that $\Phi$ induces a right action on each fiber, $\Phi_x=\Phi|_{Y_x\times\mathcal G_x}\colon Y_x\times \mathcal G_x\rightarrow Y_x$. The fibered action is said to be \emph{free} if $y\cdot g=y$ for some $(y,g)\in Y\times_X\mathcal G$ implies that $g=1_x$, $x=\pi_{Y,X}(y)$. In the same way, it is said to be \emph{proper} if the bundle morphism $Y\times_X \mathcal G\ni (y,g)\mapsto(y,y\cdot g)\in  Y\times_X Y$ is proper. If $\Phi$ is free and proper, so is each action $\Phi_x$, since the fibers of a bundle are closed. 

As the fibered action is vertical (i.e. it covers the identity $\textrm{id}_X$), we may regard the quotient space $Y/\mathcal G$ as the disjoint union of the quotients of the fibers by the induced actions, that is,
\begin{equation*}
Y/\mathcal G=\bigsqcup_{x\in X}Y_x/\mathcal G_x=\left\{[y]_{\mathcal G}=(x,[y]_{\mathcal G_x})\colon x\in X,y\in Y_x\right\}.
\end{equation*}
Obviously, the following diagram is commutative: 

\begin{equation} \label{eq:diagramafibrados}
\begin{array}{cc}
\begin{tikzpicture}
\matrix (m) [matrix of math nodes,row sep=3em,column sep=3em,minimum width=2em]
{	Y & & X \\
	& Y/\mathcal G & \\};
\path[-stealth]
(m-1-1) edge [] node [above] {$\pi_{Y,X}$} (m-1-3)
(m-1-1) edge [] node [left] {$\pi_{Y,Y/\mathcal G}\,$} (m-2-2)
(m-2-2) edge [] node [right] {$\;\,\pi_{Y/\mathcal G,X}$} (m-1-3);
\end{tikzpicture}
& 
\begin{tikzpicture}
\matrix (m) [matrix of math nodes,row sep=3em,column sep=3em,minimum width=2em]
{	y & & x \\
	& \left[y\right]_{\mathcal G} & \\};
\path[-stealth]
(m-1-1) edge [|->,decoration={markings,mark=at position 1 with {\arrow[scale=1.7]{>}}},
    postaction={decorate},shorten >=0.4pt] node [above] {} (m-1-3)
(m-1-1) edge [|->,decoration={markings,mark=at position 1 with {\arrow[scale=1.7]{>}}},
    postaction={decorate},shorten >=0.4pt] node [left] {} (m-2-2)
(m-2-2) edge [|->,decoration={markings,mark=at position 1 with {\arrow[scale=1.7]{>}}},
    postaction={decorate},shorten >=0.4pt] node [right] {} (m-1-3);
\end{tikzpicture}
\end{array}
\end{equation}

\begin{example}[Jet lift of fibered actions]\label{example:jetliftaction}
Let $\Phi\colon Y\times_X\mathcal G\to Y$ be a (right) fibered action of a Lie group bundle $\pi_{\mathcal G,X}$ on a fiber bundle $\pi_{Y,X}$. The first jet extension of $\Phi$ turns out to be a (right) fibered action of $\pi_{J^1 \mathcal G,X}$ on $\pi_{J^1 Y,X}$,
\begin{equation*}
\begin{array}{rccc}
\Phi^{(1)}\colon & J^1 Y\times_X J^1 \mathcal G & \longrightarrow & J^1 Y\\
& \left(j^1_x s,j^1_x\gamma\right) & \longmapsto & j^1_x(\Phi\circ(s,\gamma))
\end{array}
\end{equation*}
If we regard 1-jets as differentials of sections at a point, that is, $j_x^1s\equiv(ds)_x$ and $j_x^1\gamma\equiv(d\gamma)_x$ for certain local sections $s$ and $\gamma$, then $\Phi^{(1)}$ may be seen as
\begin{equation*}
\Phi^{(1)}(j^1_x s,j^1_x\gamma)\equiv(d\Phi)_{\left(s(x),\gamma(x)\right)}\circ\left((ds)_x,(d\gamma)_x\right)\colon T_x X\longrightarrow T_{s(x)\cdot\gamma(x)} Y.
\end{equation*}
\end{example}

%%%%%%%%%%%%%%%%%%
\subsection{Smooth structure of fibered quotients}

\begin{theorem}\label{theorem:Y-Y/Gfibrado}
If $\pi_{\mathcal G,X}$ acts on $\pi_{Y,X}$ freely and properly, then $Y/\mathcal G$ admits a unique smooth structure such that 
\begin{enumerate}[(i)]
    \item $\pi_{Y,Y/\mathcal G}$ is a fiber bundle with typical fiber $G$, called \emph{generalized principal bundle}.
    \item $\pi_{Y/\mathcal G,X}$ is a fibered manifold, i.e. a surjective submersion.
\end{enumerate}
\end{theorem}

\begin{proof}
Taking trivializations of $\pi_{Y,X}$ and $\pi_{\mathcal G,X}$ on the same neighborhood $\mathcal U\subset X$,  the local expression of the actions is:
\[
\mathcal U \times\hat Y\times G \longrightarrow \mathcal U \times \hat Y,
\]
where $\hat Y$ is the typical fiber of $Y$. This can be seen as a  standard Lie group action acting trivially on $\mathcal U$. The classical results on Lie group actions (for example, see \cite[Theorem 7.10]{lee2012}) provide a unique smooth structure on $(\mathcal U\times \hat Y)/G$ that can be used as a chart on $Y/\mathcal G$. In these trivializations, the projections $\pi_{Y,Y/\mathcal G}$ and $\pi_{Y/\mathcal G,X}$ are $\mathcal U\times \hat Y \to (\mathcal U \times \hat Y)/G$ and $(\mathcal U\times \hat Y)/G \to \mathcal U$ which gives \emph{(i)} and \emph{(ii)}. 
\end{proof}

\begin{remark}
If a Lie group $G$ acts freely, properly and fiberwisely on a bundle $\pi_{Y,X}\colon Y\to X$, then $\pi_{Y/G,X}\colon Y/G\to X$ is a bundle with typical fiber $\hat{Y}/G$. However, in the case of action of Lie group bundle, since the action depends on $X$, the notion of typical fiber is more delicate. For example, if $X$ is not connected, the topology of the typical fiber may differ on each component. This is the case of the action of $X\times \mathbb{Z}$ on $X\times (S^1\times \mathbb{R})$ with $X=\mathbb{R}-\{0\}$ defined as
\begin{equation*}
\begin{array}{rccc}
\Phi: & Y\times_X\mathcal G & \longrightarrow & Y\\
& \big((t,(\cos\theta,\sin\theta),z),(t,n)\big) & \longmapsto & \displaystyle \left(t,\left(\mathrm{sgn}(t)^n\cos\theta,\sin\theta\right),z+n\right).
\end{array}
\end{equation*}
The typical fiber is a Klein bottle for $t<0$ and a torus for $t>0$. 

For a connected base manifold $X$, the quotient $\hat Y/G$ is well-defined (up to diffeomorphism). Indeed, given any two points $p,q\in X$ connected by a path, $p=\gamma (0), q=\gamma(1)$, the action on $Y_t=\pi^{-1}_{Y,X}(\gamma(t))$ 
can be regarded (from the local trivializations of $\mathcal{G}$ and $Y$) as a homotopy of diffeomorphisms of $G$ on $\hat{Y}$. This gives a diffeomorphims of the fibers of $Y/\mathcal{G}$ on $p$ and $q$.
\end{remark}

We can build a trivializing atlas for $\pi_{Y,Y/\mathcal G}$ starting from a trivializing atlas $$\left\{\left(\mathcal U_\alpha,\psi_{\mathcal G}^\alpha\right)\colon\alpha\in\Lambda\right\}$$ of $\pi_{\mathcal G,X}$. For $\alpha\in\Lambda$, let $\mathcal V_\alpha=\pi^{-1}_{Y/\mathcal G,X}(\mathcal U_\alpha)$ and pick a local section $\hat s\in\Gamma(\mathcal V_\alpha,\pi_{Y,Y/\mathcal G})$ (the existence of that local section may require the choice of a smaller $\mathcal U_\alpha$). We define
\begin{equation}\label{eq:trivializacionYY/G}
\begin{array}{rccc}
\psi^\alpha_Y\colon & Y|_{\mathcal V_\alpha}=\pi^{-1}_{Y,Y/\mathcal G}(\mathcal V_\alpha) & \longrightarrow & \mathcal V_\alpha\times G\\
& y & \longmapsto & \left([y]_{\mathcal G},\hat h\right),
\end{array}
\end{equation}
where $\hat h\in G$ is such that $y=\hat s([y]_{\mathcal G})\cdot g$, with $g=(\psi^\alpha_{\mathcal G})^{-1}\left(\pi_{Y/\mathcal G,X}([y]_{\mathcal G}),\hat h\right)\in \mathcal G$. The element $\hat{h}$ exists since $Orb(y)=\pi_{Y,Y/\mathcal G}^{-1}\left([y]_{\mathcal G}\right)$ and it is unique because the fibered action is free. It turns out that $(\mathcal V_\alpha,\psi^\alpha_Y)$ is a trivialization for $\pi_{Y,Y/\mathcal G}$. Indeed, its inverse is given by
\begin{equation*}
\begin{array}{rccl}
(\psi^\alpha_Y)^{-1}\colon & \mathcal V_\alpha\times G & \longrightarrow & Y|_{\mathcal V_\alpha}\\
& \left([y]_{\mathcal G},\hat h\right) & \longmapsto & \hat s([y]_{\mathcal G})\cdot(\psi^\alpha_{\mathcal G})^{-1}\left(\pi_{Y/\mathcal G,X}([y]_{\mathcal G}),\hat h\right).
\end{array}
\end{equation*}
As a result, $\left\{(\mathcal V_\alpha,\psi_Y^\alpha)\mid \alpha\in\Lambda\right\}$ is an atlas for $\pi_{Y,Y/\mathcal G}$. 

Observe that, by definition of Lie group bundle, the maps $\psi_{\mathcal G}^\alpha|_{\mathcal G_x}\colon \mathcal G_x\rightarrow\{x\}\times G$, $x\in\mathcal U_\alpha$, are Lie group homomorphisms. It is thus straightforward that the fibered action is locally given by the right multiplication on the Lie group $G$.

\begin{corollary}\label{corollary:acciontrivializada}
Let $x\in \mathcal U_\alpha$, $y=(\psi_Y^\alpha)^{-1}\left([y]_{\mathcal G},\hat h\right)\in Y_x$ and $g=(\psi_{\mathcal G}^\alpha)^{-1}\left(x,\hat g\right)\in \mathcal G_x$, then
\begin{equation}
y\cdot g=(\psi_Y^\alpha)^{-1}\left([y]_{\mathcal G},\hat h\hat g\right).
\end{equation}
\end{corollary}

Moreover, for each $x\in\mathcal U_\alpha$ the following diagram is commutative
\begin{equation*}
\begin{tikzpicture}
\matrix (m) [matrix of math nodes,row sep=3em,column sep=5em,minimum width=2em]
{	\mathcal G_x & \{x\}\times G \\
	 \varfrak g_x & \{x\}\times\mathfrak{g} \\};
\path[-stealth]
(m-1-1) edge [] node [above] {$\psi_{\mathcal G}^\alpha|_{\mathcal G_x}$} (m-1-2)
(m-2-1) edge [] node [left] {$\exp$} (m-1-1)
(m-2-1) edge [] node [above] {$(d\psi_{\mathcal G}^\alpha|_{\mathcal G_x})_{1_x}$} (m-2-2)
(m-2-2) edge [] node [right] {$\left(\textrm{\normalfont id}_{\{x\}},\exp\right)$} (m-1-2);
\end{tikzpicture}
\end{equation*}
where $\pi_{ \varfrak g,X}$ is the Lie algebra bundle of $\pi_{\mathcal G,X}$, $\mathfrak g$ is the Lie algebra of $G$ and $\exp$ is the exponential map between the Lie algebra and its corresponding Lie group. Furthermore, note that $(d\psi_{\mathcal G}^\alpha|_{\mathcal G_x})_{1_x}$ is a linear isomorphism, since $\psi_{\mathcal G}^\alpha|_{\mathcal G_x}$ is a diffeomorphism. In other words, for each $\xi=(d\psi_{\mathcal G}^\alpha|_{\mathcal G_x})_{1_x}^{-1}\left(x,\hat\xi\right)\in\varfrak g_x$ we have
\begin{equation}\label{eq:exponencialtrivializada}
\exp\left(\xi\right)=(\psi_{\mathcal G}^\alpha)^{-1}\left(x,\exp\,\hat\xi\right).
\end{equation}
In fact, we can define a linear trivialization $(\mathcal U_\alpha,\psi_{ \varfrak g}^\alpha)$ for $\pi_{\varfrak g,X}$ from $(\mathcal U_\alpha,\psi_{\mathcal G}^\alpha)$. Namely,
\begin{equation}\label{eq:trivializacionmathfrakg}
\begin{array}{rccl}
\psi_{ \varfrak g}^\alpha\colon &  \varfrak g|_{\mathcal U_\alpha} & \longrightarrow & \mathcal U_\alpha\times\mathfrak g\\
& \xi & \longmapsto & (d\psi_{\mathcal G}^\alpha|_{\mathcal G_x})_{1_x}(\xi),\qquad x=\pi_{ \varfrak g,X}(\xi).
\end{array}
\end{equation}

%%%%%%%%%
\subsection{Infinitesimal generators}

If we fix $x\in X$, $y_0\in Y_x$ and $g_0\in \mathcal G_x$, we can consider the maps
\begin{equation*}
\begin{array}{rcclcrccl}
\Phi_{y_0}\colon & \mathcal G_x & \longrightarrow & Y_x, & \quad & \Phi_{g_0}\colon & Y_x & \longrightarrow & Y_x\\
& g & \longmapsto & y_0\cdot g, & & & y & \longmapsto & y\cdot g_0.
\end{array}
\end{equation*}
In the same way, we denote by $L_{g_0}\colon \mathcal G_x\rightarrow \mathcal G_x$ and $R_{g_0}\colon \mathcal G_x\rightarrow \mathcal G_x$ the left and right multiplication by $g_0\in \mathcal G_x$, respectively. \emph{Infinitesimal generators} (or \emph{fundamental fields}) are defined in the same fashion as in classical actions of Lie groups. Namely, for each $\xi\in\varfrak g_x$, then $\xi^*\in\mathfrak X(Y_x)$ is defined as
\begin{equation}\label{eq:definfinitesimalgenerator}
\xi^*_y=\left.\frac{d}{dt}\right|_{t=0} y\cdot\exp(t\xi)=(d\Phi _y)_{1_x}(\xi),\qquad y\in Y_x.
\end{equation}
Fundamental vector fields are $\pi_{Y,Y/\mathcal G}$-vertical, i.e. $\xi^*_y\in V_y Y=\ker(d\pi_{Y,Y/\mathcal G})_{y}$ for each $y\in Y_x$. Of course, they are also $\pi_{Y,X}$-vertical. 

\begin{proposition}\label{proposition:isomorphismverticalbundle}
Let $\pi_{\varfrak g,X}$ be the Lie algebra bundle of $\pi_{\mathcal G,X}$. The following map is a vertical isomorphism of vector bundles over $Y$: 
\begin{equation}\label{eq:isog}
\begin{array}{rcl} 
Y\times_X\varfrak g & \longrightarrow & VY \\
(y,\xi) & \longmapsto & \xi^*_y
\end{array}
\end{equation}
In addition, for any $(g,\xi)\in \mathcal G\times_X\varfrak g$, we have:
\begin{equation}
\label{eq:equivariancerighttranslation}
(\Phi _g)_*(\xi^*)=(Ad_{g^{-1}}\xi)^*.
\end{equation}
\end{proposition}
\begin{proof}
It is clear that the morphism is vertical by definition. Fix $y\in Y$ with $\pi_{Y,X}(y)=x$. Let us see that the map $Y\times_X\varfrak g_y=\varfrak g_x\ni\xi\mapsto\xi^*_y\in V_y Y$ is a linear isomorphism. The linearity is clear from the second equality of \eqref{eq:definfinitesimalgenerator}. Since $\dim{  \varfrak g_x}=\dim{\mathcal G_x}=\dim{V_y Y}$ (the last equality is for being $\pi_{Y,Y/\mathcal G}$ a submersion), we just need to prove the injectivity to conclude. Let $\xi\in\varfrak g_x$ be such that $\xi^*_y=0$. Then, $y\cdot\exp{t\xi}=y$ for every $t\in(-\epsilon,\epsilon)$. Since the action is free, this says that $\exp{t\xi}=1_x$ for every $t\in(-\epsilon,\epsilon)$ and, hence, $\xi=0$.
The second part is a straightforward computation.
\end{proof}

%%%%%%%%%
\subsection{Generalized associated bundles}

Let $\pi_{F,X}\colon F\rightarrow X$ be a fiber bundle on which $\pi_{\mathcal G,X}$ acts fiberwisely on the left. This yields a right fibered action of $\pi_{\mathcal G,X}$ on product bundle $\pi_{Y\times_X F,X}$. Namely,
\begin{equation*}
(y,f)\cdot g = (y\cdot g, g^{-1}\cdot f),\qquad (y,f,g)\in Y\times_X F\times_X\mathcal G.
\end{equation*}
If the fibered action of $\pi_{\mathcal G,X}$ on $\pi_{Y,X}$ is free and proper, so is the induced action on $\pi_{Y\times_X F,X}$. In such case, we have the smooth manifold
\begin{equation*}
Y\times _{\mathcal G} F=(Y\times_X F)/\mathcal G.
\end{equation*}
We denote the equivalence classes by $[y,f]_{\mathcal G}\in Y\times _{\mathcal G} F$.

\begin{proposition}\label{prop:fibradoasociadogeneralizado}
In the above conditions, if $\hat F$ is the typical fiber of $\pi_{F,X}$, then $\pi_{Y\times_{\mathcal G} F,Y/\mathcal G}$ is a fiber bundle with typical fiber $\hat F$, called \emph{generalized associated bundle}, where 
\begin{equation*}
\begin{array}{rccl}
\pi_{Y\times_{\mathcal G} F,Y/\mathcal G}\colon & Y\times_{\mathcal G} F & \longrightarrow & Y/\mathcal G\\
& \left[y,f\right]_{\mathcal G} &\longmapsto & [y]_{\mathcal G}
\end{array}
\end{equation*} 
\end{proposition}
\begin{proof}
It is clear that the projection is well-defined and surjective. For each $[y_0,f_0]_{\mathcal G}\in Y\times_{\mathcal G} F$, let us find a trivialization through that point. Let $\mathcal U\subset X$ be a trivializing set of both $\pi_{\mathcal G,X}$ and $\pi_{F,X}$, such that $\pi_{Y,X}(y_0)=\pi_{F,X}(f_0)\in\mathcal U$. Let $\mathcal V=\pi_{Y/\mathcal G,X}^{-1}(\mathcal U)$ and suppose that it is a trivializing set of $\pi_{Y,Y/\mathcal G}$ (maybe we need to choose a smaller $\mathcal U$ in order to achieve this). Using the above trivializations, we define
\begin{equation*}
\begin{array}{rccc}
\phi\colon & (Y\times_{\mathcal G} F)|_{\mathcal V} & \longrightarrow & \mathcal V\times\hat F\\
& [(\sigma,g),(x,f)]_{\mathcal G} & \longmapsto & \left(\sigma,g\cdot f\right),
\end{array}
\end{equation*}
where the action $G\times\hat F\to\hat F$ is induced by the fibered action of $\pi_{\mathcal G,X}$ on $\pi_{F,X}$. Observe that this action may depend on the base point $x$. Nevertheless, the condition of $\Phi\colon Y\times_X\mathcal G\to Y$ is smooth ensures that $\phi$ is also smooth. 

Of course, $[y_0,f_0]_{\mathcal G}\in (Y\times_{\mathcal G} F)|_{\mathcal V}$, since $\pi_{Y,X}(y_0)=\pi_{F,X}(f_0)\in\mathcal U$. In fact, the pair $(\mathcal V,\phi)$ is a trivialization of $\pi_{Y\times _{\mathcal G} F,Y/\mathcal G}$. Indeed, the inverse map is given by
\begin{equation*}
\begin{array}{rccc}
\phi^{-1}\colon & \mathcal V\times \hat F & \longrightarrow & (Y\times_{\mathcal G} F)|_{\mathcal V}\\
& \left(\sigma,f\right) & \longmapsto & [(\sigma,1),(x,f)]_{\mathcal G},
\end{array}
\end{equation*} 
where again $x=\pi_{Y/\mathcal G,X}(\sigma)$.
\end{proof}

\begin{example}[Conjugacy bundle]
A Lie group bundle $\pi_{\mathcal G,X}$ acts fiberwisely on the left on itself by conjugation: $g\cdot h=c_g(h)=ghg^{-1}$, $(g,h)\in\mathcal G\times_X\mathcal G$. We can thus consider the \emph{generalized conjugacy bundle}, $Y\times_{\mathcal G}\mathcal G$, which inherits the Lie group bundle structure from $\pi_{\mathcal G,X}$.
\end{example}

\begin{example}[Adjoint bundle]
In the same vein, $\pi_{\mathcal G,X}$ acts fiberwisely on the left on $\pi_{\varfrak g,X}$ via the adjoint map: $g\cdot\xi=Ad_g(\xi)=(dc_g)_{1_x}(\xi)$, $(g,\xi)\in \mathcal G\times_X\varfrak g$. The corresponding quotient is the \emph{generalized adjoint bundle}, which we will denote by $\tilde{\varfrak g}=Y\times_{\mathcal G}\varfrak g$. It is a vector bundle equipped with a Lie algebra bundle structure.
\end{example}

\begin{example}[Coadjoint bundle]
Let $\varfrak g^*$ be the dual vector bundle of $\pi_{\varfrak g,X}$, where $\pi_{\mathcal G,X}$ acts fiberwisely on the left via the coadjoint representation, $g\cdot\eta=Ad^*_{g^{-1}}(\eta)$, $(g,\eta)\in \mathcal G\times_X\varfrak g^*$. Using this action we get the \emph{generalized coadjoint bundle}, $\tilde{\varfrak g}^*=Y\times_{\mathcal G}\varfrak g^*$, which is a vector bundle.
\end{example}

%%%%%%%%%%
\section{Generalized principal connections}

%%%%%%%%%%
\subsection{Lie group bundle connections}

Recall that an Ehresmann connection (see for example \cite{michor1993}) on a fiber bundle $\pi_{Z,X}\colon Z\to X$ is a fiber map $TZ \to VZ=\ker(\pi_{Z,X})_*$ such that its restriction to $VZ$ is the identity. Similarly, we can understand an Ehresmann connection as a distribution $HZ\subset TZ$ complementary to $VZ$. Finally, an Ehresmann connection is also a section of the jet bundle $\pi_{J^1Z,Z}\colon J^1 Z \to Z$.

If $\pi_{\mathcal G,X}$ is a Lie group bundle, an Ehresmann connection
$$
\nu\colon T\mathcal G\longrightarrow V\mathcal G,\qquad U_g\longmapsto U_g^v,
$$
can be also regarded as a vertical bundle map (denoted by the same letter for the sake of simplicity)
$$\nu\colon T\mathcal G \longrightarrow \varfrak g,\qquad U_g\longmapsto \left(dR_{g^{-1}}\right)_g(U_g^v). $$
Furthermore, it is natural to impose a compatibility of $\nu$ with the algebraic structure of $\mathcal G$. 

\begin{definition}\label{def:liegroupconnection}
A \emph{Lie group bundle connection} on $\pi_{\mathcal G,X}$ is an Ehresmann connection $\nu\colon T\mathcal G \to V\mathcal G$ satisfying
\begin{enumerate}[(i)]
    \item $\ker\nu_{1_x}=(d1)_x(T_x X)$ for each $x\in X$.
    
    \item For every $(g,h)\in \mathcal G\times_X\mathcal G$ and $(U_g,U_h)\in T_g \mathcal G\times_{T_x X}T_h \mathcal G$, $x=\pi_{\mathcal G,X}(g)$, then:
\begin{equation*}
\nu_{gh}\left((dM)_{(g,h)}(U_g,U_{h})\right)=(dR_{h})_g\left(\nu_g(U_g)\right)+(dL_{g})_{h}\left(\nu_{h}(U_{h})\right),
\end{equation*}
where $M\colon\mathcal G\times_X \mathcal G\rightarrow\mathcal G$ is the fiber multiplication map.
\end{enumerate}
The corresponding conditions when regarded as a map $\nu\colon T\mathcal G\rightarrow\varfrak g$ are
\begin{enumerate}[(i)]
    \item $\ker\nu_{1_x}=(d1)_x(T_x X)$ for each $x\in X$.
    \item For every $(g,h)\in \mathcal G\times_X\mathcal G$ and $(U_g,U_h)\in T_g \mathcal G\times_{T_x X}T_h \mathcal G$, $x=\pi_{\mathcal G,X}(g)$, then:
    \begin{equation*}
    \nu_{gh}\left((dM)_{(g,h)}(U_g,U_h)\right)=\nu_g(U_g)+Ad_g\left(\nu_g(U_h)\right).
    \end{equation*}
\end{enumerate}
\end{definition}

Geometric interpretations of Lie group bundle connections are provided by the following results. We denote by $^\nu\big|\big|$ the parallel transport of $\nu$ and by $Hor_g^\nu\colon T_x X\rightarrow T_g \mathcal G$ its horizontal lift at any $g\in\mathcal G$, $x=\pi_{\mathcal G,X}(g)$.

\begin{proposition}\label{prop:liegroupconnection}
Let $\nu$ be an Ehresmann connection on $\pi_{\mathcal G,X}$ such that $\ker\nu_{1_x}=(d1)_x(T_x X)$ for each $x\in X$. Then $\nu$ is a Lie group connection if and only if for any curve $x\colon I\rightarrow X$ we have
\begin{equation}\label{eq:compatibilityofnu}
^\nu\big|\big|^{x(b)}_{x(a)}(gh)=\left( ^\nu\big|\big|^{x(b)}_{x(a)}g\right)\left( ^\nu\big|\big|^{x(b)}_{x(a)}h\right),\qquad g,h\in \mathcal G_{x(a)}.
\end{equation}
Consequently,
\begin{equation}\label{eq:inversatransportenu}
^\nu\big|\big|^{x(b)}_{x(a)}\,g^{-1}=\left(^\nu\big|\big|^{x(b)}_{x(a)}\,g\right)^{-1},\qquad^\nu\big|\big|^{x(b)}_{x(a)}\,1_{x(a)}=1_{x(b)}.
\end{equation}
\end{proposition}

\begin{proof}
Suppose that $\nu$ is a Lie group connection. Let $x\colon I\rightarrow X$ be a curve and $g,h\in \mathcal G_{x(a)}$. Denote $\alpha_1(t)={^\nu\big|}\big|^{x(t)}_{x(a)}g$, $\alpha_2(t)={^\nu\big|}\big|^{x(t)}_{x(a)}h$, $t\in I$, and $\alpha=M\circ(\alpha_1,\alpha_2)$. To show that $\alpha(t)={^\nu\big|}\big|^{x(t)}_{x(a)}(gh)$, $t\in I$, we use the uniqueness of the parallel transport. It is clear that $\pi_{\mathcal G,X}\circ\alpha=x$ and $\alpha(a)=gh$, so we only need to check that it is horizontal,
\begin{equation*}
\begin{array}{ccl}
\nu_{\alpha(t)}\left(\alpha'(t)\right) & = & \nu_{\alpha(t)}\left((dM)_{(\alpha_1(t),\alpha_2(t))}\left(\alpha_1'(t),\alpha_2'(t)\right)\right)\\
& = &\left(dR_{\alpha_2(t)}\right)_{\alpha_1(t)}\left(\nu_{\alpha_1(t)}\left(\alpha_1'(t)\right)\right)+\left(dL_{\alpha_1(t)}\right)_{\alpha_2(t)}\left(\nu_{\alpha_2(t)}\left(\alpha_2'(t)\right)\right)\;\;=\;\;0,
\end{array}
\end{equation*}
since $\alpha_1$ and $\alpha_2$ are horizontal.

Conversely, suppose that $\nu$ satisfies \eqref{eq:compatibilityofnu}. Let $x\in X$, $g,h\in \mathcal G_x$ and $(U_g,U_h)\in T_g \mathcal G\times_{T_x X}T_h \mathcal G$. Then there exists $\alpha=(\alpha_1,\alpha_2)\colon(-\epsilon,\epsilon)\rightarrow \mathcal G\times_X \mathcal G$ such that $\alpha'(0)=(\alpha_1'(0),\alpha_2'(0))=(U_g,U_h)$. Denote $x=\pi_{\mathcal G,X}\circ\alpha_1=\pi_{\mathcal G,X}\circ\alpha_2$. We conclude that $\nu$ is a Lie group connection as a straightforward consequence of \eqref{eq:compatibilityofnu} and the definition of covariant derivative,
\begin{equation*}
\begin{array}{ccl}
\nu_{gh}\left((dM)_{(g,h)}(U_g,U_{h})\right) & = & \displaystyle\nu_{gh}\left((M\circ\alpha)'(0)\right)\\
& = & \displaystyle\frac{D^\nu(M\circ\alpha)(0)}{Dt}\\
& = & \displaystyle(dR_h)_g\left(\frac{D^\nu\alpha_1(0)}{Dt}\right)+(dL_g)_h\left(\frac{D^\nu\alpha_2(0)}{Dt}\right)\\
& = & \displaystyle(dR_h)_g\left(\nu_g(U_g)\right)+(dL_g)_h\left(\nu_h(U_h)\right)
\end{array}
\end{equation*}
\end{proof}

\begin{proposition}\label{prop:liegroupconnectionjet}
Let $\nu$ be an Ehresmann connection on $\pi_{\mathcal G,X}$ and consider the corresponding jet section $\hat{\nu}\in\Gamma(\pi_{J^1 \mathcal G,\mathcal G})$. Then $\nu$ is a Lie group bundle connection if and only if
\begin{enumerate}[(i)]
    \item $\hat\nu\circ 1=j^1 1=d1$,
    \item for each $(g,h)\in \mathcal G\times_X\mathcal G$, we have that $\hat{\nu}(gh)=\hat{\nu}(g)\,\hat{\nu}(h)$ with respect to the Lie group bundle structure in $J^1\mathcal{G}$.
\end{enumerate}
\end{proposition}

\begin{proof}
It is clear that $(i)$ is equivalent to the condition: $\ker\nu_{1_x}=(d1)_x(T_x X)$ for each $x\in X$. Thanks to the previous Proposition, it is enough to show that $(ii)$ is equivalent to \eqref{eq:compatibilityofnu}. Let $x\in X$, $U_x\in T_x X$ and $g,h\in\mathcal G_x$ and pick a curve $\alpha\colon(-\epsilon,\epsilon)\rightarrow X$ such that $\alpha'(0)=U_x$. We define the following curve:
\begin{equation*}
\beta_g(t)={^\nu\big|\big|}_{\alpha(0)}^{\alpha(t)}g,\qquad t\in(-\epsilon,\epsilon).
\end{equation*}
We have $\beta_g'(0)=Hor_g^\nu(U_x)=\hat\nu(g)(U_x)$, where $Hor_g^\nu\colon T_x X\rightarrow T_g \mathcal G$ is the horizontal lifting given by $\nu$. Performing the same construction with $h$ and $gh$ we obtain curves $\beta_h,\beta_{gh}\colon(-\epsilon,\epsilon)\rightarrow T\mathcal G$.

Suppose that \eqref{eq:compatibilityofnu} is satisfied, then $\beta_{gh}=M\circ(\beta_g,\beta_h)$. Therefore
\begin{equation*}
\begin{array}{ccl}
\hat\nu(gh)(U_x) & = & \beta_{gh}'(0)\\
& = & (dM)_{(g,h)}\left(\beta_g'(0),\beta_h'(0)\right)\\
& = & (dM)_{(g,h)}\left(\hat\nu(g)(U_x),\hat\nu(h)(U_x)\right)
\end{array}
\end{equation*}
Since the equality is valid for every $U_x\in T_x X$, we conclude that $\hat\nu(gh)=\hat\nu(g)\,\hat\nu(h)$.

Conversely, if $(ii)$ holds, then $\beta_{gh}(0)=gh=\beta_g(0)\,\beta_h(0)$ and $\beta_{gh}'(0)=(dM)_{(g,h)}(\beta_g'(0),\beta_h'(0))=\left(M\circ(\beta_g,\beta_h)\right)'(0)$. By uniqueness of the parallel transport we obtain that $\beta_{gh}=M\circ(\beta_g,\beta_h)$ for each local curve $\alpha\colon(-\epsilon,\epsilon)\rightarrow X$ such that $\alpha(0)=x$.
\end{proof}

From the proof above we also deduce that the condition
\begin{equation*}
Hor_{gh}^\nu=(dM)_{(g,h)}\circ\left(Hor_g^\nu,Hor_h^\nu\right),\qquad (g,h)\in\mathcal G\times_X\mathcal G,
\end{equation*}
together with $(i)$ of Definition \ref{def:liegroupconnection} characterize Lie group connections. More generally, we have the following property.

\begin{proposition}
Let $\nu$ be a Lie group connection on $\pi_{\mathcal G,X}$. Then for each $(g,h)\in \mathcal G\times_X\mathcal G$, $U_x\in T_x X$, $x=\pi_{\mathcal G,X}(g)$, and $U_h\in T_h \mathcal G$ such that $(d\pi_{\mathcal G,X})_h(U_h)=U_x$ we have
\begin{equation*}
(dM)_{(g,h)}\left(Hor_g^\nu(U_x),U_h\right)=Hor_{gh}^\nu(U_x)+\left(dL_g\right)_h\left(\nu_h(U_h)\right).
\end{equation*}
\end{proposition}

\begin{proof}
Thanks to the uniqueness of the horizontal lifting, it is enough to show that 
\begin{equation*}
U=(dM)_{(g,h)}\left(Hor_g^\nu(U_x),U_h\right)-\left(dL_g\right)_h\left(\nu_h(U_h)\right)
\end{equation*} is horizontal and projects to $U_x$. For horizontality we have that
\begin{equation*}
\begin{array}{ccl}
\nu_{gh}(U) & = & \left(dR_h\right)_g\left(\nu_g\left(Hor_g^\nu(U_x)\right)\right)+\left(dL_g\right)_h\left(\nu_h\left(U_h\right)\right)-\nu_{gh}\left(\left(dL_g\right)_h\left(\nu_h(U_h)\right)\right)\\
& = & \left(dL_g\right)_h\left(\nu_h\left(U_h\right)\right)-\nu_{gh}\left(\left(dL_g\right)_h\left(\nu_h(U_h)\right)\right)\\
& = & 0.
\end{array}
\end{equation*}
The last equality is a particular case of property $(ii)$ of the definition of Lie group connection with $U_g=0$. In such case, we have $(dM)_{(g,h)}\left(0,U_h\right)=\left(dL_g\right)\left(U_h\right)$ and, thus,
\begin{equation*}
\nu_{gh}\left(\left(dL_g\right)\left(U_h\right)\right)=\nu_{gh}\left((dM)_{(g,h)}\left(0,U_h\right)\right)=\left(dL_g\right)_h\left(\nu_h\left(U_h\right)\right).
\end{equation*}
At last, we check that it projects to $U_x$. In order to do so, let $\alpha=\left(\alpha_1,\alpha_2\right)\colon(-\epsilon,\epsilon)\rightarrow \mathcal G\times_X \mathcal G$ such that $\alpha(0)=(g,h)$ and $\alpha'(0)=\left(Hor_g^\nu(U_x),U_h\right)$,
\begin{equation*}
\begin{array}{ccl}
(d\pi_{\mathcal G,X})_{gh}(U) & = & (d\pi_{\mathcal G,X})_{gh}\left((dM)_{(g,h)}\left(Hor_g^\nu(U_x),U_h\right)\right)-(d\pi_{\mathcal G,X})_{gh}\left(\left(dL_g\right)_h\left(\nu_h(U_h)\right)\right)\vspace{0.1cm}\\
& = & \displaystyle d(\pi_{\mathcal G,X}\circ M)_{(g,h)}\left(Hor_g^\nu(U_x),U_h\right)\vspace{0.1cm}\\
& = & \displaystyle\left.\frac{d}{dt}\right|_{t=0}(\pi_{\mathcal G,X}\circ M\circ \alpha)(t)\vspace{0.1cm}\\
& = & \displaystyle\left.\frac{d}{dt}\right|_{t=0}(\pi_{\mathcal G,X}\circ\alpha_1)(t)\vspace{0.1cm}\\
& = & \displaystyle(d\pi_{\mathcal G,X})_g\left(Hor_g^\nu(U_x)\right)\vspace{0.1cm}\\
& =& U_x.
\end{array}
\end{equation*}
\end{proof}
With respect to this last Proposition, if the vector $U_h$ is horizontal, then 
\[
(dM)_{(g,h)}\left(H_g,U_h\right)=H_{gh},\qquad U_h\in T_h\mathcal{G}, \quad \pi_{\mathcal{G},X}(g)=\pi_{\mathcal{G},X}(h),
\]
an identity that is similar to the corresponding property of connections on standard principal bundles.
%%%%%%

\subsection{Induced connection on the Lie algebra bundle}

A Lie group connection induces a linear connection on the corresponding Lie algebra bundle.

\begin{proposition}\label{prop:nablaginducida}
Let $x\colon I\rightarrow X$ be a smooth curve. Then the map
${^{\varfrak g}\big|\big|}_{x(a)}^{x(b)}\colon\varfrak g_{x(a)}\rightarrow\varfrak g_{x(b)}$ defined as
\begin{equation*}
{^{\varfrak g}\big|\big|}_{x(a)}^{x(b)}\xi=\left.\frac{d}{d\epsilon}\right|_{\epsilon=0}{^\nu\big|\big|}_{x(a)}^{x(b)}\exp(\epsilon\,\xi),\quad \xi\in\varfrak g_{x(a)}
\end{equation*}
is a linear parallel transport on $\pi_{\varfrak g,X}$.
\end{proposition}

\begin{proof}
The map is well-defined, since $^\nu\big|\big|^{x(b)}_{x(a)}\,1_{x(a)}=1_{x(b)}$. Likewise, it is easy to check that ${^{\varfrak g}\big|\big|}_{x(a)}^{x(b)}\lambda\xi=\lambda{^{\varfrak g}\big|\big|}_{x(a)}^{x(b)}\xi$ for each $\xi\in\varfrak g_{x(a)}$ and $\lambda\in\mathbb R$. To conclude, let $\xi,\eta\in\varfrak g_{x(a)}$ and $\epsilon\in\mathbb R$ ``small enough''. We apply the Zassenhaus formula (for example, cf. \cite{casas2012}):
\begin{equation*}
\exp(\epsilon(\xi+\eta))=\exp(\epsilon\,\xi)\exp(\epsilon\,\eta)\prod_{n=2}^{\infty}\exp(C_n(\epsilon\,\xi,\epsilon\,\eta)),
\end{equation*}
where $C_n(\epsilon\,\xi,\epsilon\,\eta)$ is a homogeneous Lie polynomial in $\{\epsilon\,\xi,\epsilon\,\eta\}$ for each $n\geq 2$. Thus, $\exp(\epsilon(\xi+\eta))=\exp(\epsilon\,\xi)\exp(\epsilon\,\eta)\exp\left(o(\epsilon^2)\right)$ and are done:
\begin{equation*}
\begin{array}{ccl}
\displaystyle{^{\varfrak g}\big|\big|}_{x(a)}^{x(b)}(\xi+\eta) & = & \displaystyle\left.\frac{d}{d\epsilon}\right|_{\epsilon=0}{^\nu\big|\big|}_{x(a)}^{x(b)}\exp(\epsilon(\xi+\eta))\vspace{0.1cm}\\
& = & \displaystyle\left.\frac{d}{d\epsilon}\right|_{\epsilon=0}{^\nu\big|\big|}_{x(a)}^{x(b)}\left(\exp(\epsilon\,\xi)\exp(\epsilon\,\eta)\exp\left(o\left(\epsilon^2\right)\right)\right)\vspace{0.1cm}\\
& = & \displaystyle\left.\frac{d}{d\epsilon}\right|_{\epsilon=0}\left({^\nu\big|\big|}_{x(a)}^{x(b)}\exp(\epsilon\,\xi)\right)\left({^\nu\big|\big|}_{x(a)}^{x(b)}\exp(\epsilon\,\eta)\right)\left({^\nu\big|\big|}_{x(a)}^{x(b)}\exp\left(o\left(\epsilon^2\right)\right)\right)\vspace{0.1cm}\\
& = & \displaystyle\left.\frac{d}{d\epsilon}\right|_{\epsilon=0}{^\nu\big|\big|}_{x(a)}^{x(b)}\exp(\epsilon\,\xi)+\left.\frac{d}{d\epsilon}\right|_{\epsilon=0}{^\nu\big|\big|}_{x(a)}^{x(b)}\exp(\epsilon\,\eta)+\underset{0}{\underbrace{\left.\frac{d}{d\epsilon}\right|_{\epsilon=0}{^\nu\big|\big|}_{x(a)}^{x(b)}\exp\left(o\left(\epsilon^2\right)\right)}}\vspace{0.1cm}\\
& = & \displaystyle{^{\varfrak g}\big|\big|}_{x(a)}^{x(b)}\xi+{^{\varfrak g}\big|\big|}_{x(a)}^{x(b)}\eta .
\end{array}
\end{equation*}
\end{proof}

This linear parallel transport naturally respects the adjoint representation of $\varfrak g$.

\begin{proposition}\label{3.5}
Let $x\colon I\rightarrow X$. Then
\begin{equation*}
{^{\varfrak g}\big|\big|}_{x(a)}^{x(b)}Ad_g(\xi)= Ad_{{^\nu||}_{x(a)}^{x(b)}g}\left({^{\varfrak g}\big|\big|}_{x(a)}^{x(b)}\xi\right),\qquad g\in\mathcal G_{x(a)},\quad \xi\in\varfrak g_{x(a)}.
\end{equation*}
\end{proposition}

\begin{proof}
Consider the curves $\alpha(t)=Ad_{{^\nu||}_{x(a)}^{x(t)}g}\left({^{\varfrak g}\big|\big|}_{x(a)}^{x(t)}\xi\right)$ and $\beta(t)={^{\varfrak g}\big|\big|}_{x(a)}^{x(t)}Ad_g(\xi)$, $t\in I$, which satisfy $\alpha(a)=\beta(a)=Ad_g(\xi)$ and $\pi_{\varfrak g,X}\circ\alpha=\pi_{\varfrak g,X}\circ\beta=x$. Thanks to the uniqueness of parallel transport, to show that $\alpha=\beta$ it is enough to check that $\alpha'=\beta'$. For each $t\in I$, using the definition of $^{\varfrak g}\big|\big|$ we have
\begin{equation*}
\begin{array}{ccl}
\alpha'(t) %& = & \displaystyle\frac{d}{dt}\;Ad_{{^\nu||}_{x(a)}^{x(t)}g}\left({^{\varfrak g}\big|\big|}_{x(a)}^{x(t)}\xi\right)\vspace{0.1cm}\\
%& = & \displaystyle\frac{d}{dt}\left(dc_{{^\nu||}_{x(a)}^{x(t)}g}\right)_{1_x}\left(\left.\frac{d}{d\epsilon}\right|_{\epsilon=0}{^\nu\big|\big|}_{x(a)}^{x(t)}\exp(\epsilon\,\xi)\right)\vspace{0.1cm}\\
& = & \displaystyle\frac{d}{dt}\left.\frac{d}{d\epsilon}\right|_{\epsilon=0}c_{{^\nu||}_{x(a)}^{x(t)}g}\left({^\nu\big|\big|}_{x(a)}^{x(t)}\exp(\epsilon\,\xi)\right)\vspace{0.1cm}\\
& \overset{(\star)}{=} & \displaystyle\frac{d}{dt}\left.\frac{d}{d\epsilon}\right|_{\epsilon=0}{^\nu\big|\big|}_{x(a)}^{x(t)}c_g\left(\exp(\epsilon\,\xi)\right)\vspace{0.1cm}\\
& \overset{(\star\star)}{=} & \displaystyle\frac{d}{dt}\left.\frac{d}{d\epsilon}\right|_{\epsilon=0}{^\nu\big|\big|}_{x(a)}^{x(t)}\exp\left((dc_g)_{1_{x(a)}}(\epsilon\,\xi)\right)\vspace{0.1cm}\\
%& = & \displaystyle\frac{d}{dt}\left.\frac{d}{d\epsilon}\right|_{\epsilon=0}{^\nu\big|\big|}_{x(a)}^{x(t)}\exp\left(\epsilon Ad_g(\xi)\right)\vspace{0.1cm}\\
%& = & \displaystyle\frac{d}{dt}\;{^{\varfrak g}\big|\big|}_{x(a)}^{x(t)}Ad_g(\xi)\vspace{0.1cm}\\
& = & \displaystyle\beta'(t),
\end{array}
\end{equation*}
where $(\star)$ is due to $\nu$ being a Lie group connection and $(\star\star)$ comes from the fact that $c_g\colon \mathcal G_{x(a)}\rightarrow \mathcal G_{x(a)}$ is a Lie group homomorphism.
\end{proof}

We denote by $\nabla^{\varfrak g}$ and $\nabla^{\varfrak g}/dt$ the linear connection and the covariant derivative on $\pi_{\varfrak g,X}$ corresponding to this parallel transport ${^{\varfrak g}\big|\big|}$, respectively. A more practical way of regarding $\nabla^{\varfrak g}$ may be the following. Consider the Lie group bundle connection as a map $\nu :T\mathcal{G}\to \varfrak{g}$. Then, for any section $\xi\in\Gamma(\pi _{\varfrak g,X})$, we have
\begin{equation*}
\nabla^{\varfrak g}\xi =\left.\frac{d}{dt}\right|_{t=0} \nu\circ d\,\mathrm{exp}(t\xi).
\end{equation*}
The next result is the infinitesimal version of Proposition \ref{3.5}.

\begin{proposition}
Let $g\colon I\rightarrow \mathcal G$ and $\xi\colon I\rightarrow\varfrak g$ be such that $\pi_{\mathcal G,X}\circ g=\pi_{\varfrak g,X}\circ\xi=x$. Then for all $t\in I$ we have
\begin{equation*}
\frac{\nabla^{\varfrak g}\left(Ad_g\circ\xi\right)(t)}{dt}=Ad_{g(t)}\left(\frac{\nabla^{\varfrak g}\xi(t)}{dt}\right)+\left[\left(dR_{g(t)^{-1}}\right)_{g(t)}\left(\frac{\nabla^\nu g(t)}{dt}\right),Ad_{g(t)}(\xi(t))\right]_{\mathcal{G}}.
\end{equation*}
\end{proposition}

\begin{proof}
Fix $t\in I$ and denote $x=\pi_{\mathcal G,X}\circ g$. Let $\alpha\colon(-\epsilon,\epsilon)\rightarrow \mathcal G_{x(t)}$ be the curve given by $\alpha(s)={^\nu\big|\big|}_{x(t+s)}^{x(t)} g(t+s)$ for each $s\in(-\epsilon,\epsilon)$. We have $\alpha(0)=g(t)$ and $\alpha'(0)=\nabla^\nu g(t)/dt$. Denoting $\beta(s)=\alpha(s)\,g(t)^{-1}=R_{g(t)^{-1}}\left(\alpha(s)\right)$, $s\in(-\epsilon,\epsilon)$, we have that $\beta(0)=1_{x(t)}$ and $\beta'(0)=\left(dR_{g(t)^{-1}}\right)_{g(t)}\left(\nabla^\nu g(t)/dt\right)$. Hence,
\begin{equation*}
\begin{array}{ccl}
\displaystyle\left.\frac{d}{ds}\right|_{s=0} Ad_{\alpha(s)}\left(\xi(t)\right) & = & \displaystyle\left.\frac{d}{ds}\right|_{s=0} Ad_{\beta(s)}\left(Ad_{g(t)}(\xi(t))\right)\vspace{0.1cm}\\
& = & \displaystyle ad(\beta'(0))\left(Ad_{g(t)}(\xi(t))\right)\vspace{0.1cm}\\
& = & \displaystyle\left[\left(dR_{g(t)^{-1}}\right)_{g(t)}\left(\frac{\nabla^\nu g(t)}{dt}\right),Ad_{g(t)}(\xi(t))\right]_{\mathcal{G}}
\end{array}
\end{equation*}

Using this, we conclude
\begin{equation*}
\begin{array}{ccl}
\displaystyle\frac{\nabla^{\varfrak g}\left(Ad_{g}\circ\xi\right)(t)}{dt} & = & \displaystyle\left.\frac{d}{ds}\right|_{s=0}{^{\nabla^{\varfrak g}}\big|\big|}_{x(t+s)}^{x(t)}Ad_{g(t+s)}(\xi(t+s))\vspace{0.1cm}\\
& = & \displaystyle\left.\frac{d}{ds}\right|_{s=0}Ad_{{^\nu||_{x(t+s)}^{x(t)}}g(t+s)}\left({^{\nabla^{\varfrak g}}\big|\big|}_{x(t+s)}^{x(t)}\xi(t+s)\right)\vspace{0.1cm}\\
& = & \displaystyle\left.\frac{d}{ds}\right|_{s=0}Ad_{g(t)}\left({^{\nabla^{\varfrak g}}\big|\big|}_{x(t+s)}^{x(t)}\xi(t+s)\right)+\left.\frac{d}{ds}\right|_{s=0}Ad_{{^\nu||_{x(t+s)}^{x(t)}}g(t+s)}\left(\xi(t)\right)\vspace{0.1cm}\\
& = & \displaystyle Ad_{g(t)}\left(\frac{\nabla^{\varfrak g}\xi(t)}{dt}\right)+\left[\left(dR_{g(t)^{-1}}\right)_{g(t)}\left(\frac{\nabla^\nu g(t)}{dt}\right),Ad_{g(t)}(\xi(t))\right]_{\mathcal{G}}
\end{array}
\end{equation*}
\end{proof}

%%%%%%%%%%
\subsection{Generalized principal connections}

Let $\pi_{Y,X}\colon Y\to X$ be a fiber bundle on which a Lie group bundle $\pi_{\mathcal G,X}\colon\mathcal G \to X$ acts freely and properly on the right, and denote by $\Phi\colon Y\times_X \mathcal G \to Y$ the fibered action. 

\begin{lemma}\label{lemma:transformacionverticales}
For every $(y,g)\in Y\times_X \mathcal G$ and $\xi,\eta\in\varfrak g_x$, $x=\pi_{\mathcal{G},X}(g)$, we have
\begin{equation*}
    (d\Phi)_{(y,g)}\left(\xi^*_y,\eta ^* _g\right)=(Ad_{g^{-1}}\left(\xi+\eta \right))^*_{y\cdot g},
\end{equation*}
where $\eta_g^*=d/dt|_{t=0}\exp(t\eta)g$ is the infinitesimal generator of $\eta$ at $g$.
\end{lemma}

\begin{proof}
It is a straightforward computation using that
\begin{equation*}
Ad_{g^{-1}}\left(\eta_g\right)^*_{y\cdot g}=\left.\frac{d}{dt}\right|_{t=0}\Phi _y\left(\exp{(t\,\eta_g)}\,g\right).
\end{equation*}

\begin{comment}
Observe that
\begin{equation*}
\begin{array}{ccl}
Ad_{g^{-1}}\left(\nu(U_g)\right)^*_{y\cdot g} & = & \left.\frac{d}{dt}\right|_{t=0}y\cdot g\,\exp{(t\,Ad_{g^{-1}}(\nu(U_g)))}\vspace{0.1cm}\\
& = & \left.\frac{d}{dt}\right|_{t=0}y\cdot g\cdot\exp{\left(t\, \left.\frac{d}{ds}\right|_{s=0}c_{g^{-1}}(\exp(s\,\nu(U_g)))\right)}\vspace{0.1cm}\\
& = & \left.\frac{d}{dt}\right|_{t=0}y\cdot g\, c_{g^{-1}}\left(\exp{(t\,\nu(U_g))}\right)\vspace{0.1cm}\\
& = & \left.\frac{d}{dt}\right|_{t=0}y\cdot \exp{(t\,\nu(U_g))}\, g\vspace{0.1cm}\\
& = & \left.\frac{d}{dt}\right|_{t=0}\Phi _y\left(\exp{(t\,\nu(U_g))}\,g\right)
\end{array}
\end{equation*}
Subsequently,
\begin{equation*}
\begin{array}{ccl}
(d\Phi)_{(y,g)}(X^*_y,\nu_g(U_g)) & = & \left.\frac{d}{dt}\right|_{t=0}\Phi\left(y\cdot\exp{(tX)},\exp{(t\nu(U_g)g)}\right)\vspace{0.1cm}\\
& = & \left.\frac{d}{dt}\right|_{t=0}\Phi\left(y\cdot\exp{(tX)},g\right)+\left.\frac{d}{dt}\right|_{t=0}\Phi\left(y,\exp{(t\nu(U_g))}\, g\right)\vspace{0.1cm}\\
& = & \left.\frac{d}{dt}\right|_{t=0}\Phi _g\left(y\cdot\exp(tX)\right)+\left.\frac{d}{dt}\right|_{t=0}\Phi _y\left(\exp{(t\nu(U_g))}\,g\right)\vspace{0.1cm}\\
& = & (d\Phi _g)_y(X^*_y)+Ad_{g^{-1}}\left(\nu(U_g)\right)^*_{y\cdot g}\vspace{0.1cm}\\
& = & Ad_{g^{-1}}\left(X\right)^*_{y\cdot g}+Ad_{g^{-1}}\left(\nu(U_g)\right)^*_{y\cdot g}\vspace{0.1cm}\\
& = & Ad_{g^{-1}}\left(X+\nu(U_g)\right)^*_{y\cdot g}.
\end{array}
\end{equation*}
\end{comment}
\end{proof}

Ehresmann connections on $\pi_{Y,Y/\mathcal G}$ are identified with forms as follows.

\begin{proposition}\label{prop:equivalenciaconexiones}
There is a bijective correspondence between Ehresmann connections on $\pi_{Y,Y/\mathcal G}$ and $1$-forms\footnote{In fact, we are considering forms with values in 
$(\pi_{Y,X})^*(\varfrak g)=Y\times_X\varfrak g$ but, for the sake of simplicity, we do not write it. In the same fashion, abusing the notation we suppose that $\omega_y\in T_y^*Y\otimes\varfrak g_x$ for each $y\in Y$, $x=\pi_{Y,X}(y)$.}
$\omega\in\Omega^1(Y,\varfrak g)$ such that
\begin{equation*}
\omega_y(\xi^*_y)=\xi,\qquad (y,\xi)\in Y\times_X\varfrak g.
\end{equation*}
Furthermore,  such forms satisfy
\begin{equation*}
\omega_{y\cdot g}\left((d\Phi _g)_y(\xi_y^*)\right)=Ad_{g^{-1}}\left(\omega_y(\xi_y^*)\right),\qquad (y,g,\xi)\in Y\times_X\mathcal G\times_X\varfrak g.    
\end{equation*}
\end{proposition}

\begin{proof}
The equivalence between Ehresmann connections on $\pi_{Y,Y/\mathcal G}$ and the 1-forms as in the statement is a consequence of the isomorphism \eqref{eq:isog}. 

For the second part, let $(y,g,\xi)\in Y\times_X\mathcal G\times_X\varfrak g$, then:
\begin{equation*}
\begin{array}{l}
\omega_{y\cdot g}\left((d\Phi _g)_y(\xi_y^*)\right)=\omega_{y\cdot g}\left((Ad_{g^{-1}}\xi)_{y\cdot g}^*\right)=Ad_{g^{-1}}(\xi)=Ad_{g^{-1}}\left(\omega_y(\xi_y^*)\right),
\end{array}
\end{equation*}
where equation \eqref{eq:equivariancerighttranslation} has been taken into account.
\end{proof}

Observe that $U_y^v=\omega_y(U_y)^*_y$ for each $U_y\in T_y Y$, $y\in Y$, and, thus $U_y$ is horizontal if and only if $\omega_y(U_y)=0$. We are ready to introduce generalized principal connections.

\begin{definition}\label{def:connectiondistribution}
Let $\nu$ be an Ehresmann connection on the Lie group bundle $\pi_{\mathcal G,X}$. A \emph{generalized principal connection} on $\pi_{Y,Y/\mathcal G}$ associated to $\nu$ is an Ehresmann connection $H\subset TY$ on $\pi_{Y,Y/\mathcal G}$ satisfying:
\begin{equation*}
\left[(d\Phi)_{(y,g)}\left(U_y,U_g\right)\right]^v=(d\Phi _g)_y(U_y^v)+(dL_{g^{-1}})_g(\nu_g(U_g))^*_{y\cdot g},
\end{equation*}
for every $(y,g)\in Y\times_X \mathcal G$ and $(U_y,U_g)\in T_yY\times_{T_x X} T_g \mathcal G$, where $x=\pi_{Y,X}(y)$.
\end{definition}

Thanks to Proposition \ref{prop:equivalenciaconexiones}, it is easy to check that the following is an equivalent way of defining a generalized principal connection.

\begin{definition}\label{def:connectionform}
Let $\nu$ be an Ehresmann connection on the Lie group bundle $\pi_{\mathcal G,X}$. A \emph{generalized principal connection} on $\pi_{Y,Y/\mathcal G}$ associated to $\nu$ is a form $\omega\in\Omega^1(Y,\varfrak g)$ satisfying:

\begin{enumerate}[(i)]
    \item (Complementarity) $\omega_y(\xi^*_y)=\xi$ for every $(y,\xi)\in Y\times_X\varfrak g$.
    \item ($Ad$-equivariance) For each $(y,g)\in Y\times_X \mathcal G$ and $(U_y,U_g)\in T_yY\times_{T_x X} T_g \mathcal G$, $x=\pi_{Y,X}(x)$, then
\begin{equation*}
\omega_{y\cdot g}\left((d\Phi)_{(y,g)}(U_y,U_g)\right)=Ad_{g^{-1}}\left(\omega_y(U_y)+\nu_g(U_g)\right).
\end{equation*}
\end{enumerate}
\end{definition}

Roughly speaking, generalized principal connections extend the property given in Lemma \ref{lemma:transformacionverticales} and Proposition \ref{prop:equivalenciaconexiones} to non-necessarily vertical vectors $U_g$ with respect to $\nu$.

The next result gives a geometric interpretation of the above Definitions in terms of the parallel transports ${^\nu\big|}\big|$ and ${^\omega\big|}\big|$, in the same vein as Proposition \ref{prop:liegroupconnection}.

\begin{proposition}\label{prop:compatibilityofomega}
Let $\nu\colon T\mathcal G\to\varfrak g$ and $\omega\in\Omega^1(Y,\varfrak g)$ be Ehresmann connections on $\pi_{\mathcal G,X}$ and $\pi_{Y,Y/\mathcal G}$, respectively. Then $\omega$ is a generalized principal connection associated to $\nu$ if and only if for any curve $\gamma\colon I\rightarrow Y/\mathcal G$, the corresponding parallel transports satisfy
\begin{equation}\label{eq:compatibilityofomega}
{^\omega\big|}\big|^{\gamma(t)}_{\gamma(a)}(y\cdot g)=\left({^\omega\big|}\big|^{\gamma(t)}_{\gamma(a)}y\right)\cdot\left( {^\nu\big|}\big|^{x(t)}_{x(a)}g\right),\qquad g\in \mathcal G_{x(a)},\quad y\in Y_{\gamma(a)},\quad t\in I,
\end{equation}
where $x=\pi_{Y/\mathcal G,X}\circ\gamma$.
\end{proposition}

\begin{proof}
Suppose that $\omega$ is a generalized principal connection associated to $\nu$. Let $\gamma\colon I\rightarrow Y/\mathcal G$ be a curve, $x=\pi_{Y/\mathcal G,X}\circ\gamma$, $y\in Y_{\gamma(a)}$ and $g\in \mathcal G_{x(a)}$. Denote $\alpha_1(\cdot)={^\omega\big|}\big|^{\gamma(\cdot)}_{\gamma(a)}y$, $\alpha_2(\cdot)={^\nu\big|}\big|^{x(\cdot)}_{x(a)}g$ and $\alpha=\Phi\circ(\alpha_1,\alpha_2)$. To show that $\alpha(\cdot)={^\omega\big|}\big|^{\gamma(\cdot)}_{\gamma(a)}(y\cdot g)$ we use the uniqueness of the parallel transport. It is clear that $\pi_{Y,Y/\mathcal G}\circ\alpha=\gamma$ and $\alpha(a)=y\cdot g$, so we only need to check that $\alpha$ is horizontal. For each $t\in I$ we have
\begin{equation*}
\begin{array}{ccl}
\alpha'(t)^v & = & (d\Phi)_{(\alpha_1(t),\alpha_2(t))}\left(\alpha_1'(t),\alpha_2'(t)\right)^v\vspace{0.1cm}\\
& = &\left(d\Phi_{\alpha_2(t)}\right)_{\alpha_1(t)}\left(\alpha_1'(t)^v\right)+\left(dL_{\alpha_2(t)^{-1}}\right)_{\alpha_2(t)}\left(\nu_{\alpha_2(t)}\left(\alpha_2'(t)\right)\right)^*_{\alpha(t)}\;\;=\;\;0
\end{array}
\end{equation*}
since $\alpha_1$ and $\alpha_2$ are horizontal for the corresponding connections.

Conversely, suppose that $\omega$ and $\nu$ satisfy \eqref{eq:compatibilityofomega}. Let $(y,g)\in Y\times_X \mathcal G$ and $(U_y,U_g)\in T_y Y\times_{T_x X}T_g$. We consider $\alpha=(\alpha_1,\alpha_2)\colon(-\epsilon,\epsilon)\rightarrow Y\times_X \mathcal G$ such that $\alpha'(0)=(\alpha_1'(0),\alpha_2'(0))=(U_y,U_g)$ and we denote $\gamma=\pi_{Y,Y/\mathcal G}\circ\alpha_1$ and $x=\pi_{Y,X}\circ\alpha_1=\pi_{\mathcal G,X}\circ\alpha_2$. A straightforward consequence of \eqref{eq:compatibilityofomega} and the definition of covariant derivative is the following:
\begin{equation*}
\begin{array}{ccl}
(\Phi\circ\alpha)'(0)^v & = &  \displaystyle\left.\frac{D^\omega(\Phi\circ\alpha)}{Dt}\right|_{t=0}\vspace{0.1cm}\\
& = & \displaystyle(d\Phi_g )_y\left(\left.\frac{D^\omega\alpha_1}{Dt}\right|_{t=0}\right)+(d\Phi_y )_g\left(\left.\frac{D^\nu\alpha_2}{Dt}\right|_{t=0}\right)\vspace{0.1cm}\\
& = & \displaystyle(d\Phi_g )_y\left(U_y^v\right)+(d\Phi_y )_g\left(U_g^v\right)
\end{array}
\end{equation*}
Furthermore, note that $(d\Phi_y )_g\left(U_g^v\right)=(d\Phi_{y\cdot g} )_{1_x}\left((dL_{g^{-1}})_g\left(U_g^v\right)\right)=(dL_{g^{-1}})_g\left(U_g^v\right)^*_{y\cdot g}$. Thus, we conclude that $\omega$ is a generalized principal connection associated to $\nu$:
\begin{equation*}
(d\Phi)_{(y,g)}\left(U_y,U_g\right)^v=(\Phi\circ\alpha)'(0)^v=(d\Phi_g )_y\left(U_y^v\right)+(dL_{g^{-1}})_g\left(U_g^v\right)^*_{y\cdot g}.
\end{equation*}
\end{proof}

\begin{proposition}\label{prop:jetsectionomega}
Let $\nu\colon T\mathcal G\to\varfrak g$ and $\omega\in\Omega^1(Y,\varfrak g)$ be Ehresmann connections on $\pi_{\mathcal G,X}$ and $\pi_{Y,Y/\mathcal G}$, respectively, and suppose that $Y/\mathcal G=X$. Let $\hat\omega\in\Gamma(\pi_{J^1 Y,Y})$ and $\hat\nu\in\Gamma(\pi_{J^1\mathcal G,G})$ be the corresponding jet sections. Then $\omega$ is a generalized principal connection associated to $\nu$ if and only if 
\begin{equation}\label{eq:equivariancejetsections}
\hat\omega(y\cdot g)=\hat\omega(y)\cdot\hat\nu(g),\qquad(y,g)\in Y\times_X\varfrak g,
\end{equation}
where the fibered action of $\pi_{J^1\mathcal G,X}$ on $\pi_{J^1 Y,X}$ is given by the first jet extension of $\Phi$ (recall Example \ref{example:jetliftaction}).
\end{proposition}

\begin{proof}
Thanks to Proposition \ref{prop:compatibilityofomega} we only need to show that \eqref{eq:equivariancejetsections} is equivalent to \eqref{eq:compatibilityofomega}. Let $(y,g)\in Y\times_X\mathcal G$, $x=\pi_{Y,X}(y)$ and $U_x\in T_x X$, and pick a curve $\alpha\colon(-\epsilon,\epsilon)\to X$ with $\alpha'(0)=U_x$. We define the curves
\begin{equation*}
\beta_y(t)={^\omega\big|\big|}_{\alpha(0)}^{\alpha(t)}y,\quad\beta_{y\cdot g}(t)={^\omega\big|\big|}_{\alpha(0)}^{\alpha(t)}(y\cdot g),\quad\beta_y(t)={^\nu\big|\big|}_{\alpha(0)}^{\alpha(t)}g,
\end{equation*}
This way, we have $\beta_g'(0)=Hor_g^{\nu}(U_x)=\hat\nu(g)(U_x)$ and $\beta_y'(0)=Hor_y^{\omega}(U_x)=\hat\omega(y)(U_x)$, and analogous for $\beta_{y\cdot g}$.

Suppose that \eqref{eq:compatibilityofomega} is satisfied. Then $\beta_{y\cdot g}=\Phi\circ(\beta_y,\beta_g)$. Hence
\begin{equation*}
\begin{array}{ccl}
\hat\omega(y\cdot g)(U_x) & = & \beta_{y\cdot g}'(0)\\
& = & (d\Phi)_{(y,g)}\left(\beta_y'(0),\beta_g'(0)\right)\\

& = & (d\Phi)_{(y,g)}\left(\hat\omega(y)(U_x),\hat\nu(g)(U_x)\right).
\end{array}
\end{equation*}
This equation is valid for every $U_x\in T_x X$, whence $\hat\omega(y\cdot g)=\hat\omega(y)\cdot\hat\nu(g)$.

Conversely, if \eqref{eq:equivariancejetsections} is satisfied, then $\beta_{y\cdot g}(0)=y\cdot g=\beta_y(0)\cdot\beta_g(0)=(\Phi\circ(\beta_y,\beta_g))(0)$ and
\begin{equation*}
\begin{array}{ccl}
\beta_{y\cdot g}'(0) & = & \hat\omega(y\cdot g)(U_x)\\
& = &\left(\hat\omega(y)\cdot\hat\nu(g)\right)(U_x)\\
& = & (d\Phi)_{(y,g)}\left(\hat\omega(y)(U_x),\hat\nu(g)(U_x)\right)\\
& = & (d\Phi)_{(y,g)}\left(\beta_y'(0),\beta_g'(0)\right)\\
& = & (\Phi\circ(\beta_y,\beta_g))'(0).
\end{array}
\end{equation*}
By uniqueness of the parallel transport we conclude that $\beta_{y\cdot g}=\Phi\circ(\beta_y,\beta_g)$.
\end{proof}
\color{black}

We denote by $Hor^\omega_{y}\colon T_{[y]_{\mathcal G}}(Y/\mathcal G)\rightarrow T_y Y$ the \emph{horizontal lift} given by $\omega$ at $y\in Y$. A different geometric interpretation of generalized principal connections is the following.

\begin{proposition}\label{prop:transformacionlevantamientohorizontal}
Let $(y,g)\in Y\times_X \mathcal G$ and $(U_{[y]_{\mathcal G}},U_g)\in T_{[y]_{\mathcal G}}(Y/\mathcal G)\times_{T_x X}T_g\mathcal G$, where $x=\pi_{Y,X}(y)$. Then 
\begin{equation*}
(d\Phi)_{(y,g)}\left(Hor^\omega_y(U_{[y]_{\mathcal G}}),U_g\right)=Hor^\omega_{y\cdot g}(U_{[y]_{\mathcal G}})+(dL _{g^{-1}})_g(\nu_g(U_g))^*_{y\cdot g}
\end{equation*}
\end{proposition}

\begin{proof}
We just need to prove that $U=(d\Phi)_{(y,g)}\left(Hor^\omega_y(U_{[y]_{\mathcal G}}),U_g\right)-(dL _{g^{-1}})_g(\nu_g(U_g))^*_{y\cdot g}$ is horizontal and projects to $U_{[y]_{\mathcal G}}$. With respect to the former
\begin{equation*}
\begin{array}{ccl}
U^v & = & (d\Phi)_{(y,g)}\left(Hor^\omega_y(U_{[y]_{\mathcal G}}),U_g\right)^v-(dL_{g^{-1}})_g(\nu_g(U_g))^*_{y\cdot g}\vspace{0.1cm}\\
& = & (d\Phi _g)_y\left(Hor^\omega_y(U_{[y]_{\mathcal G}})^v\right)\vspace{0.1cm}\\
& = & 0
\end{array}
\end{equation*}
For the latter let  $\alpha=(\alpha_1,\alpha_2)\colon(-\epsilon,\epsilon)\rightarrow Y\times_X \mathcal G$ be such $\alpha'(0)=\left(Hor^\omega_y(U_{[y]_{\mathcal G}}),U_g\right)$. We then have
\begin{equation*}
\begin{array}{ccl}
(d\pi_{Y,Y/\mathcal G})_{y\cdot g}(U) & = & \displaystyle(d\pi_{Y,Y/\mathcal G})_{y\cdot g}\left((d\Phi)_{(y,g)}\left(Hor^\omega_y(U_{[y]_{\mathcal G}}),U_g\right)\right)\vspace{0.1cm}\\
& = & \displaystyle d(\pi_{Y,Y/\mathcal G}\circ\Phi)_{(y,g)}\left(Hor^\omega_y(U_{[y]_{\mathcal G}}),U_g\right)\vspace{0.1cm}\\
& = & \displaystyle\left.\frac{d}{dt}\right|_{t=0}(\pi_{Y,Y/\mathcal G}\circ\Phi\circ\alpha)(t)\vspace{0.1cm}\\
& = & \displaystyle\left.\frac{d}{dt}\right|_{t=0}(\pi_{Y,Y/\mathcal G}\circ\alpha_1)(t)\vspace{0.1cm}\\
& = & \displaystyle(d\pi_{Y,Y/\mathcal G})_y\left(Hor^\omega_y(U_{[y]_{\mathcal G}})\right)\vspace{0.1cm}\\
& = & U_{[y]_{\mathcal G}}
\end{array}
\end{equation*}
\end{proof}

\begin{theorem}[Existence of generalized principal connections]\label{theorem:existencia}
If $X$ is a paracompact smooth manifold, then there exist an Ehresmann connection $\nu$ on $\pi_{\mathcal G,X}$ and a generalized principal connection $\omega$ on $\pi_{Y,Y/\mathcal G}$ associated to $\nu$.
\end{theorem}

\begin{proof}
Let $\{(\mathcal U_\alpha,\psi_{\mathcal G}^\alpha)\mid\alpha\in\Lambda\}$ be a trivializing atlas for $\pi_{\mathcal G,X}$, $\{(\mathcal V_\alpha=\pi_{Y/\mathcal G,X}^{-1}(\mathcal U_\alpha),\psi_Y^\alpha)\mid\alpha\in\Lambda\}$ be the induced trivializing atlas for $\pi_{Y,Y/\mathcal G}$, as in \eqref{eq:trivializacionYY/G}, and $\{(\mathcal U_\alpha,\psi_{\varfrak g}^\alpha)\mid\alpha\in\Lambda\}$ be the induced trivializing atlas for $\pi_{\varfrak g,X}$, as in \eqref{eq:trivializacionmathfrakg}. For $\hat g\in G$, we denote by $\hat L_{\hat g}\colon G\rightarrow G$ and $\hat R_{\hat g}\colon G\to G$ the left and right multiplication by $\hat g$, respectively. In the same way, we denote by $\widehat{Ad}$ the adjoint representation of $G$.

Fixed $\alpha\in\Lambda$, we define the local Ehresmann connection $\nu_\alpha\colon T\mathcal G|_{\mathcal U_\alpha}\rightarrow V\mathcal G|_{\mathcal U_\alpha}$ on $\pi_{\mathcal G,X}|_{\mathcal U_\alpha}\colon\mathcal G|_{\mathcal U_\alpha}\rightarrow \mathcal U_\alpha$ as follows
\begin{equation*}
(\nu_\alpha)_g\left((d\psi_{\mathcal G}^\alpha)_g^{-1}\left(U_x,\left(d\hat L_{\hat g}\right)_1\left(\hat\eta\right)\right)\right)=(d\psi_{\mathcal G}^\alpha)_g^{-1}\left(0_x,\left(d\hat L_{\hat g}\right)_1\left(\hat\eta\right)\right),\qquad U_x\in T_x X,\quad\hat\eta\in \mathfrak g
\end{equation*}
where $x\in \mathcal U_\alpha$, $\hat g\in G$ and $g=(\psi_{\mathcal G}^\alpha)^{-1}\left(x,\hat g\right)\in\mathcal G|_{\mathcal U_\alpha}$. Likewise, we define the local 1-form $\omega_\alpha\in\Omega^1\left(Y|_{\mathcal V_\alpha},\varfrak g|_{\mathcal U_\alpha}\right)$ as
\begin{equation*}
(\omega_\alpha)_y\left((d\psi_Y^\alpha)_y^{-1}\left(U_{[y]_{\mathcal G}},\left(d\hat L_{\hat h}\right)_1\left(\hat\xi\right)\right)\right)=(\psi_{\varfrak g}^\alpha)^{-1}\left(x,\hat\xi\right),\qquad U_{[y]_{\mathcal G}}\in T_{[y]_{\mathcal G}}(Y/\mathcal G),\quad\hat\xi\in \mathfrak g,
\end{equation*}
where $[y]_{\mathcal G}\in \mathcal V_\alpha$, $\hat h\in G$ and $y=(\psi_Y^\alpha)^{-1}\left([y]_{\mathcal G},\hat h\right)\in Y|_{\mathcal V_\alpha}$. We show that $\nu_\alpha$ and $\omega_\alpha$ satisfy both properties of Definition \ref{def:connectionform}:

\begin{enumerate}[(i)]
    \item (Complementarity) Let $\xi=(\psi_{\varfrak g}^\alpha)^{-1}\left(x,\hat\xi\right)\in\varfrak g|_{\mathcal U_\alpha}$ and $y=(\psi_Y^\alpha)^{-1}\left([y]_{\mathcal G},\hat h\right)\in Y|_{\mathcal V_\alpha}$. Using Corollary \ref{corollary:acciontrivializada} and equation (\ref{eq:exponencialtrivializada}) we get
    \begin{equation}\label{eq:generadorinfinitesimaltrivializado}
    \xi^*_y = %\left.\frac{d}{dt}\right|_{t=0}\Phi(y,\exp(t\xi)) = \left.\frac{d}{dt}\right|_{t=0}(\psi^\alpha_Y)^{-1}\left([y]_{\mathcal G},\hat h\exp(t\hat\xi)\right) = 
    (d\psi_Y^\alpha)_y^{-1}\left(0_{[y]_{\mathcal G}},\left(d\hat L_{\hat h}\right)_{1}\left(\hat\xi\right)\right)
    \end{equation}
    Thence,
    \begin{equation*}
    (\omega_\alpha)_y\left(\xi_y^*\right)=(\omega_\alpha)_y\left((d\psi_Y^\alpha)_y^{-1}\left(0_{[y]_{\mathcal G}},\left(d\hat L_{\hat h}\right)_{1}\left(\hat\xi\right)\right)\right)=(\psi_{\varfrak g}^\alpha)^{-1}\left(x,\hat\xi\right)=\xi
    \end{equation*}

    \item (Ad-equivariance) Let $\hat h,\hat g\in G$ and $[y]_{\mathcal G}\in \mathcal V_\alpha$, and consider $y=(\psi_Y^\alpha)^{-1}\left([y]_{\mathcal G},\hat h\right)$ and $g=(\psi_{\mathcal G}^\alpha)^{-1}\left(x,\hat g\right)$, where $x=\pi_{Y/\mathcal G,X}\left([y]_{\mathcal G}\right)\in \mathcal U_\alpha$. Likewise, let $\hat\xi,\hat\eta\in\mathfrak g$ and $U_{[y]_{\mathcal G}}\in T_{[y]_{\mathcal G}}(Y/\mathcal G)$, and consider
    \begin{equation*}
    U_y=(d\psi_Y^\alpha)_y^{-1}\left(U_{[y]_{\mathcal G}},\left(d\hat L_{\hat h}\right)_1\left(\hat\xi\right)\right)\in T_y Y,\qquad U_g=(d\psi_{\mathcal G}^\alpha)_g^{-1}\left(U_x,\left(d\hat L_{\hat g}\right)_1\left(\hat\eta\right)\right)\in T_g \mathcal G,    
    \end{equation*}
    where $U_x=(d\pi_{Y/\mathcal G,X})_{[y]_{\mathcal G}}\left(U_{[y]_{\mathcal G}}\right)\in T_x X$. Consider curves $\gamma\colon(-\epsilon,\epsilon)\rightarrow Y|_{\mathcal V_\alpha}$ and $\beta\colon(-\epsilon,\epsilon)\rightarrow \mathcal G|_{\mathcal U_\alpha}$ such that $\gamma'(0)=U_y$ and $\beta'(0)=U_g$. Then we have that $\psi_Y^\alpha\circ\gamma=(\gamma_1,\gamma_2)$ and $\psi_{\mathcal G}^\alpha\circ\beta=(\beta_1,\beta_2)$ for certain curves $\gamma_1\colon(-\epsilon,\epsilon)\rightarrow \mathcal V_\alpha$, $\gamma_2\colon(-\epsilon,\epsilon)\rightarrow G$, $\beta_1\colon(-\epsilon,\epsilon)\rightarrow \mathcal U_\alpha$ and $\beta_2\colon(-\epsilon,\epsilon)\rightarrow G$ such that $\gamma_1'(0)=U_{[y]_{\mathcal G}}$, $\gamma_2'(0)=\left(d\hat L_{\hat h}\right)_1\left(\hat\xi\right)$, $\beta_1'(0)=U_x$ and $\beta_2'(0)=\left(d\hat L_{\hat g}\right)_1\left(\hat\eta\right)$. Note that we can always choose the curves satisfying $\pi_{Y/\mathcal G,X}\circ\gamma_1=\beta_1$. Using these curves, Corollary \ref{corollary:acciontrivializada} and Equation \eqref{eq:generadorinfinitesimaltrivializado} it can be seen that
    \begin{equation*}
    \begin{array}{ccl}
    (d\Phi)_{(y,g)}\left(U_y,U_g\right) & = & \displaystyle\left.\frac{d}{dt}\right|_{t=0}\Phi(\alpha(t),\beta(t))\\
    %& = & \left.\frac{d}{dt}\right|_{t=0}(\psi_Y^\alpha)^{-1}\left(\gamma_1(t),\gamma_2(t)\beta_2(t)\right)\\
    %& = & (d\psi_Y^\alpha)_{y\cdot g}^{-1}\left(U_{[y]_{\mathcal G}},(d\hat R_{\hat g})_{\hat h}\left((d\hat L_{\hat h})_1(\hat\xi)\right)+(d\hat L_{\hat h})_{\hat g}\left((d\hat L_{\hat g})_1(\hat\eta)\right)\right)\\
    %& = & (d\psi_Y^\alpha)_{y\cdot g}^{-1}\left(U_{[y]_{\mathcal G}},(d\hat R_{\hat g})_{\hat h}\left((d\hat L_{\hat h})_1(\hat\xi)\right)\right)\\
    %& & +(d\psi_Y^\alpha)_{y\cdot g}^{-1}\left(0_{[y]_{\mathcal G}},(d\hat L_{\hat h})_{\hat g}\left((d\hat L_{\hat g})_1(\hat\eta)\right)\right)\\
    %& = & (d\psi_Y^\alpha)_{y\cdot g}^{-1}\left(U_{[y]_{\mathcal G}},(d\hat L_{\hat h\hat g})_{1}\left(\widehat{Ad}_{\hat g^{-1}}(\hat\xi)\right)\right)\\
    %& & +(d\psi_Y^\alpha)_{y\cdot g}^{-1}\left(0_{[y]_{\mathcal G}},(d\hat L_{\hat h\hat g})_{1}(\hat\eta)\right)\\
    & = & \displaystyle(d\psi_Y^\alpha)_{y\cdot g}^{-1}\left(U_{[y]_{\mathcal G}},\left(d\hat L_{\hat h\hat g}\right)_{1}\left(\widehat{Ad}_{\hat g^{-1}}\left(\hat\xi\right)\right)\right) +(\psi_{\varfrak g}^\alpha)^{-1}\left(x,\hat\eta\right)^*_{y\cdot g}
    \end{array}
    \end{equation*}
    This, together with property (i), lead to
    \begin{equation*}
    \begin{array}{ccl}
    (\omega_\alpha)_{y\cdot g}\left((d\Phi)_{y\cdot g}(U_y,U_g)\right) & = & (\psi_{\varfrak g}^\alpha)^{-1}\left(x,\widehat{Ad}_{\hat g^{-1}}\left(\hat\xi\right)\right)+(\psi_{\varfrak g}^\alpha)^{-1}\left(x,\hat\eta\right)\\
    & = & Ad_{g^{-1}}\left((\psi_{\varfrak g}^\alpha)^{-1}\left(x,\hat\xi\right)\right)+(d\psi_{\mathcal G}^\alpha)_{1_x}^{-1}\left(0_x,\hat\eta\right)\\
    %& = & Ad_{g^{-1}}\left((\omega_\alpha)_y\left((d\psi_Y^\alpha)_y^{-1}\left(U_{[y]_{\mathcal G}},(d\hat L_{\hat h})_1(\hat\xi)\right)\right)\right)\\
    %& & +(dL_{g^{-1}})_g\left((d\psi_{\mathcal G}^\alpha)_{g}^{-1}\left(0_x,(d\hat L_{\hat g})_1(\hat\eta)\right)\right)\\
     & = & Ad_{g^{-1}}\left((\omega_\alpha)_y(U_y)\right)+Ad_{g^{-1}}\left(\left(dR_{g^{-1}}\right)_g\left((\nu_\alpha)_g(U_g)\right)\right)\\
     & = & Ad_{g^{-1}}\left((\omega_\alpha)_y(U_y)+(\nu_\alpha)_g(U_g)\right),
    \end{array}
    \end{equation*}
    where we have used \eqref{eq:trivializacionmathfrakg} and the fact that $(d\psi_{\mathcal G}^\alpha)_{1_x}^{-1}\left(0_x,\hat\eta\right)=(d\psi_{\mathcal G}^\alpha|_{\mathcal G_x})_{1_x}^{-1}(x,\hat\eta)$.
\end{enumerate}

Therefore, $\omega_\alpha$ is a generalized principal connection on $\pi_{Y,Y/\mathcal G}|_{\mathcal V_\alpha}\colon Y|_{\mathcal V_\alpha}\rightarrow \mathcal V_\alpha$ associated to $\nu_\alpha$. At last, we take a smooth partition of unity $\{\theta_\alpha\mid\alpha\in\Lambda\}$ on $X$ subordinated to $\{\mathcal U_\alpha\mid\alpha\in\Lambda\}$. It is easy to check that $\{\Theta_\alpha=\theta_\alpha\circ\pi_{Y/\mathcal G,X}\mid\alpha\in\Lambda\}$ is a smooth partition of unity on $Y/\mathcal G$ subordinated to $\{\mathcal V_\alpha\mid\alpha\in\Lambda\}$. 
Now set $\omega=\sum_{\alpha\in\Lambda}(\Theta_\alpha\circ\pi_{Y,Y/\mathcal G})\,\omega_\alpha$ and $\nu=\sum_{\alpha\in\Lambda}(\theta_\alpha\circ\pi_{\mathcal G,X})\,\nu_\alpha$. They are a well defined 1-form $\omega\in\Omega^1(Y,\varfrak g)$ and a well defined Ehresmann connection $\nu\colon T\mathcal G\rightarrow V\mathcal G$. It is straightforward that $\nu$ and $\omega$ satisfy both properties of Definition \ref{def:connectionform}. Thus, $\omega$ is a generalized principal connection associated to $\nu$.
\end{proof}

\begin{comment}
\begin{definition}
A form $\alpha\in\Omega^p(Y,\varfrak g)$ is said to be \textbf{tensorial of the adjoint type} if it is

\begin{enumerate}[(i)]
    \item \textbf{Horizontal}, i.e. $\alpha_y(U_1,...,U_p)=0$ whenever one $U_i\in V_y Y=\ker(d\pi_{Y,Y/\mathcal G})_y$, for every $y\in Y$ and $U_1,...,U_p\in T_y Y$.
    
    \item \textbf{$Ad$-equivariant}, i.e. for every $(y,g)\in Y\times_X \mathcal G$, $(U^1_y,U^1_g),...,(U^p_y,U^p_g)\in T_y Y\times_{T_x X}T_g \mathcal G$, $x=\pi_{Y,X}(y)$, then
\begin{equation*}
\alpha_{y\cdot g}\left((d\Phi)_{(y,g)}(U^1_y,U^1_g),...,(d\Phi)_{(y,g)}(U^p_y,U^p_g)\right)=Ad_{g^{-1}}\left(\alpha_y(U^1_y,...,U^p_y)\right)
\end{equation*} 
\end{enumerate}
\textcolor{red}{Creo que no es necesario tomar los $U_g^i$ horizontales. Por una parte, el concepto de horizontalidad no está definido en principio, porque no hemos tenido que fijar necesariamente una conexión generalizada. En cualquier caso, suponiendo que lo hemos hecho, siempre podemos descomponer $U_g^i=(U_g^i)^h+(U_g^i)^v$, y así
\begin{equation*}
(d\Phi)_{(y,g)}\left(U_y^i,U_g^i\right)=(d\Phi)_{(y,g)}\left(U_y^i,(U_g^i)^h\right)+(d\Phi)_{(y,g)}\left(0,(U_g^i)^v\right)
\end{equation*}
El segundo sumando es un vector vertical de $T_{y\cdot g} Y$, por la propia definición de conexión generalizada. Por lo tanto, al aplicarle $\alpha$ se anula, por ser $\alpha$ una forma horizontal.}
\end{definition}
\end{comment}

A 1-form $\alpha\in\Omega^1(Y,\varfrak g)$ is said to be \emph{tensorial of the adjoint type} if it is \emph{horizontal}, i.e. $\alpha_y(U_y)=0$ for each $U_y\in V_y Y=\ker(d\pi_{Y,Y/\mathcal G})_y$, $y\in Y$, and
\emph{$Ad$-equivariant}, i.e.
\begin{equation*}
\alpha_{y\cdot g}\left((d\Phi)_{(y,g)}(U_y,U_g)\right)=Ad_{g^{-1}}\left(\alpha_y(U_y)\right)
\end{equation*}
for each $(U_y,U_g)\in T_y Y\times_{T_x X}T_g \mathcal G$, $(y,g)\in Y\times_X \mathcal G$, $x=\pi_{Y,X}(y)$. The family of tensorial 1-forms of the adjoint type is denoted by $\overline\Omega^1(Y,\varfrak g)$. These forms can be reduced to $Y/\mathcal G$ in the sense that  we have a bijection 
\begin{equation*}
\begin{array}{ccc}
\overline\Omega^1(Y,\varfrak g) & \longrightarrow & \Omega^1(Y/\mathcal G,\tilde{\varfrak g})\\
\alpha & \longmapsto & \tilde\alpha,
\end{array}   
\end{equation*}
given by $\tilde\alpha_{[y]_{\mathcal G}}(U_{[y]_{\mathcal G}})=[y,\alpha_y(U_y)]_{\mathcal G}$, where $U_y\in T_y Y$ projects to $U_{[y]_{\mathcal G}}\in T_{[y]_{\mathcal G}}(Y/\mathcal G)$, $y\in Y$. By an abuse of notation, we identify $\alpha\equiv\tilde\alpha$. This construction can be straightforwardly generalized to $p$-forms.

\begin{proposition}\label{prop:espacioafinconexiones}
Let $\nu\colon T\mathcal G\rightarrow\varfrak g$ be an Ehresmann connection on $\pi_{\mathcal G,X}$. The family of generalized principal connections on $\pi_{Y,Y/\mathcal G}$ associated to $\nu$ is an affine space modelled on $\overline\Omega^1(Y,\varfrak g)$.
\end{proposition}

\begin{proof}
On the one hand, if $\omega_1,\omega_2\in\Omega^1(Y,\varfrak g)$ are generalized principal connections on $\pi_{Y,Y/\mathcal G}$ associated to $\nu$, then $\omega_1-\omega_2\in\overline\Omega^1(Y,\varfrak g)$. The $Ad$-equivariance comes from the $Ad$-equivariance of $\omega_1$ and $\omega_2$ and the fact that the adjoint map is linear. For the horizontality, let $y\in Y$ and $V_y\in V_y Y=\ker{(d\pi_{Y,Y/\mathcal G})_y}$, we have
\begin{equation*}
(\omega_1-\omega_2)_y(V_y)^*_y=(\omega_1)_y(V_y)^*_y-(\omega_2)_y(V_y)^*_y=V_y-V_y=0.
\end{equation*}
As the map \eqref{eq:isog} is an isomorphism, we deduce that $(\omega_1-\omega_2)_y(V_y)=0$. 

On the other hand, given a generalized principal connection $\omega\in\Omega^1(Y,\varfrak g)$ associated to $\nu$ and a 1-form $\overline\omega\in\overline\Omega^1(Y,\varfrak g)$, a straightforward computation shows that $\omega+\overline\omega$ is a generalized principal connection, since $\overline\omega$ vanish on vertical vectors.
\end{proof}

Now we give another interpretation of a generalized principal connection as a connection on $\pi_{Y\times_X\mathcal G,X}$ equivariant under the fibered action.

\begin{proposition}
Let $\omega\in\Omega^1(Y,\varfrak g)$ be a generalized principal connection on $\pi_{Y,Y/\mathcal G}$ associated to an Ehresmann connection $\nu\colon T\mathcal G\rightarrow V\mathcal G$ on $\pi_{\mathcal G,X}$. Then the map 
\begin{equation*}
\begin{array}{rccl}
\varpi\colon & T(Y\times_X \mathcal G) & \longrightarrow & V(Y\times_X \mathcal G)\\
& (U_y,U_g) & \longmapsto & \left(\omega_y(U_y)_y^*,\nu_g(U_g)\right),
\end{array}
\end{equation*}
where $V(Y\times_X \mathcal G)=\ker\left(\pi_{Y\times_X \mathcal G,X}\right)_*$, is an Ehresmann connection on $\pi_{Y\times_X \mathcal G,X}$ equivariant under the map $\mathbf\Phi\colon Y\times_X \mathcal G\rightarrow Y\times_X \mathcal G$ defined as $\mathbf\Phi(y,g)=(y\cdot g,g)$, i.e.
\begin{equation*}
    (d\mathbf\Phi)_{(y,g)}\left(\varpi_{(y,g)}\left(U_y,U_g\right)\right)=\varpi_{\mathbf\Phi(y,g)}\left((d\mathbf\Phi)_{y,g}(U_y,U_g)\right), 
\end{equation*}
for each $(U_y,U_g)\in T_y Y\times_{T_x X}T_g \mathcal G,\quad (y,g)\in Y\times_X \mathcal G$.
\end{proposition}

\begin{proof}
To begin with, observe that $\varpi$ is well defined, i.e. $\varpi_{(y,g)}(U_y,U_g)\in V_{(y,g)}(Y\times_X \mathcal G)$ for every $(U_y,U_g)\in T_y Y\times_{T_x X} T_g \mathcal G$, $(y,g)\in Y\times_X \mathcal G$. This is a straightforward consequence of equality $V_{(y,g)}(Y\times_X \mathcal G)=V_y Y\times V_g\mathcal G$, where $V_y Y=\ker(\pi_{Y,X})_*$, and the fact that infinitesimal generators are $\pi_{Y,X}$-vertical. Likewise, it is clear that $\varpi$ is a vertical vector bundle morphism over $Y\times_X \mathcal G$, since $\omega_y$, $\nu_g $ and the infinitesimal generator are linear maps. To conclude, let us check the $\mathbf\Phi$-equivariance. Note that $(d\mathbf\Phi)_{(y,g)}\left(U_y,U_g\right)=\left((d\Phi)_{(y,g)}(U_y,U_g),U_g\right)$. Hence,
\begin{equation*}
\begin{array}{ccl}
(d\mathbf\Phi)_{(y,g)}\left(\varpi_{(y,g)}\left(U_y,U_g\right)\right) & = & (d\mathbf\Phi)_{(y,g)}\left(\omega_y(U_y)_y^*,\nu_g(U_g)\right)\vspace{0.1cm}\\
& = & \left((d\Phi)_{(y,g)}\left(\omega_y(U_y)_y^*,\nu_g(U_g)\right),\nu_g(U_g)\right)\vspace{0.1cm}\\
& = & \left(Ad_{g^{-1}}\left(\omega_y(U_y)+\nu_g(U_g)\right)^*_{\Phi(y,g)},\nu_g(U_g)\right)\vspace{0.1cm}\\
& = & \left(\omega_{\Phi(y,g)}\left((d\Phi)_{(y,g)}(U_y,U_g)\right)^*_{\Phi(y,g)},\nu_g(U_g)\right)\vspace{0.1cm}\\
& = & \varpi_{\left(\Phi(y,g),g\right)}\left((d\Phi)_{(y,g)}\left(U_y,U_g\right),U_g\right)\vspace{0.1cm}\\
& = & \varpi_{\mathbf\Phi(y,g)}\left((d\mathbf\Phi)_{(y,g)}\left(U_y,U_g\right)\right)
\end{array}
\end{equation*}
\end{proof}

So far, we have considered generalized principal connections associated to arbitrary Ehresmann connections on $\mathcal{G}\to X$. We now prove that this Ehresmann connections must be a Lie group bundle connection, that is, they must respect the algebraic structure of $\pi_{\mathcal G,X}$.

\begin{proposition}
Let $\omega\in\Omega^1(Y,\varfrak g)$ be a generalized principal connection on $\pi_{Y,Y/\mathcal G}$ associated to an Ehresmann connection $\nu\colon T\mathcal G\to\varfrak g$ on $\pi_{\mathcal G,X}$. Then $\nu$ is a Lie group bundle connection.
\end{proposition}

\begin{proof}
Let $\gamma\colon I\to Y/\mathcal G$ projecting onto $x\colon I\to X$. Thanks to \eqref{eq:compatibilityofomega} for each $y\in Y_{\gamma(a)}$ and $g,h\in\mathcal G_{x(a)}$ we have
\begin{equation*}
\begin{array}{ccl}
\left({^\omega\big|}\big|^{\gamma(b)}_{\gamma(a)}y\right)\cdot\left( {^\nu\big|}\big|^{x(b)}_{x(a)}g\right)\left({^\nu\big|}\big|^{x(b)}_{x(a)}h\right) & = & \left({^\omega\big|}\big|^{\gamma(b)}_{\gamma(a)}(y\cdot g)\right)\cdot\left({^\nu\big|}\big|^{x(b)}_{x(a)}h\right)\vspace{0.1cm}\\
& = & {^\omega\big|}\big|^{\gamma(b)}_{\gamma(a)}(y\cdot gh)\vspace{0.1cm}\\
& = & \left({^\omega\big|}\big|^{\gamma(b)}_{\gamma(a)}y\right)\cdot\left({^\nu\big|}\big|^{x(b)}_{x(a)}(gh)\right).
\end{array}
\end{equation*}
Since the action is free, we conclude that property \eqref{eq:compatibilityofnu} holds.

On the other hand, applying \eqref{eq:compatibilityofomega} with $g=1_{x(a)}$ and using again the fact that the action is free, we get
\begin{equation*}
{^\nu\big|}\big|^{x(t)}_{x(a)}1_{x(a)}=1_{x(b)}.
\end{equation*}
This gives that $\ker\nu_{1_x}=(d1)_x(T_x X)$ for each $x\in X$ and we conclude thanks to Proposition \ref{prop:liegroupconnection}. Indeed, let $U_x\in\ker\nu_{1_x}=H_{1_x}$ and denote $u_x=(d\pi_{\mathcal G,X})_{1_x}(U_x)\in T_x X$. Let $x\colon(-\epsilon,\epsilon)\to X$ be such that $x'(0)=u_x$. Then we have that
\begin{equation*}
U_x=\left.\frac{d}{dt}\right|_{t=0}{^\nu\big|}\big|^{x(t)}_{x(0)}1_x=\left.\frac{d}{dt}\right|_{t=0} 1_{x(t)}=(d1)_x(u_x).
\end{equation*}
Therefore, $\ker\nu_{1_x}\subset(d1)_x(T_x X)$. The other inclusion is due to the fact that both are vector spaces of the same dimension.
\end{proof}

%%%%%%%%%%%%%%%%
\subsection{Curvature}

The \emph{cuvature} $\omega$ as an Ehresmann connection (see, for example \cite[\S 9.4]{michor1993}) is the $VY$-valued 2-form $\Omega\in\Omega^2(Y,VY)$ defined as
\begin{equation*}
\Omega(U_1,U_2)=-\left[U_1-\omega(U_1)^*,U_2-\omega(U_2)^*\right],\qquad U_1,U_2\in\mathfrak X(Y).
\end{equation*}
This is equivalent, by the identification \eqref{eq:isog}, to the $\varfrak{g}$-valued 2-form $\Omega\in\Omega^2(Y,\varfrak{g})$ 
\begin{equation*}
\Omega(U_1,U_2)=-\omega\left(\left[U_1-\omega(U_1)^*,U_2-\omega(U_2)^*\right]\right),
\end{equation*}
which will be denoted the same.

The linear connection $\nabla^{\varfrak g}$ induced on $\pi_{\varfrak g,X}$ by $\nu$ enables us to express the curvature as follows.

\begin{proposition}\label{prop:curvature}
Let $d^{\varfrak g}$ be the exterior covariant derivative\footnote{
The \emph{exterior covariant derivative} of a linear connection $\nabla^E$ on a vector bundle $\pi_{E,X}$ is an operator in the family of $E$-valued forms on $X$, $d^E\colon\Omega^\bullet(X,E)\longrightarrow \Omega^{\bullet+1}(X,E)$. For a 1-form $\alpha\in\Omega^1(X,E)$ it is given by
\begin{equation*}
d^E \alpha (U_1,U_2)=\nabla_{U_1}^E(\alpha(U_2))-\nabla_{U_2}^E(\alpha(U_1))-\alpha([U_1,U_2]),\qquad U_1,U_2\in\mathfrak X(X).
\end{equation*}
}
associated to $\nabla^{\varfrak g}$. Then\footnote{
The generalized connection $\omega$ takes values on the vector bundle $\pi_{Y\times_X\varfrak g,Y}$, which is the pull-back of $\pi_{\varfrak g,X}$ by $\pi_{Y,X}$. Hence, we also need to pull-back the linear connection $\nabla^{\varfrak g}$ to compute the exterior covariant derivative of $\omega$. Abusing the notation, we denote the pull-back connection by the same symbol: $\pi_{Y,X}^*\left(\nabla^{\varfrak g}\right)\equiv\nabla^{\varfrak g}$.
}
\begin{equation*}
\Omega\left(U_1,U_2\right)=d^{\varfrak g}\,\omega\left(U_1^h,U_2^h\right),\qquad U_1,U_2\in\mathfrak X(Y).
\end{equation*}
\end{proposition} 

As in the case of (standard) principal connections, it is possible to regard the curvature as a 2-form on the base space $Y/\mathcal G$ with values in $\tilde{\varfrak g}$. 

\begin{definition}\label{def:reducedcurvature}
The \emph{reduced curvature} of $\omega$ is the 2-form $\tilde\Omega\in\Omega^2\left(Y/\mathcal G,\tilde{\varfrak g}\right)$ taking values in the adjoint bundle given by
\begin{equation*}
\tilde\Omega_{[y]_{\mathcal G}}\left(U_1,U_2\right)=\left[y,\Omega_y\left(Hor_y^\omega(U_1),Hor_y^\omega(U_2)\right)\right]_{\mathcal G}
\end{equation*}
for each $[y]_{\mathcal G}\in Y/\mathcal G$ and $U_1,U_2\in T_{[y]_{\mathcal G}}(Y/\mathcal G)$, where $y\in Y$ is such that $\pi_{Y,Y/\mathcal G}(y)=[y]_{\mathcal G}$.
\end{definition}

The reduced curvature is well-defined, i.e. it does not depend on the choice of $y\in Y$. Indeed, let $g\in\mathcal G_x$, where $x=\pi_{Y,X}(y)$, $u_i=(d\pi_{Y/\mathcal G,X})_{[y]_{\mathcal G}}(U_i)\in T_x X$ for $i=1,2$ and $\gamma\in\Gamma(\pi_{\mathcal G,X})$ be such that $\gamma(x)=g$ and $\nu_{\gamma(x)}\circ(d\gamma)_x=0$. Proposition \ref{prop:transformacionlevantamientohorizontal} gives
\begin{equation*}
Hor_{y\cdot g}^\omega(U_i)=(d\Phi)_{(y,g)}\left(Hor_y^\omega(U_i),(d\gamma)_x(u_i)\right),\qquad i=1,2.
\end{equation*}
Hence, we have
\begin{equation}\label{eq:curvaturareducida}
\begin{array}{l}
\left[y\cdot g,\Omega_{y\cdot g}\left(Hor_{y\cdot g}^\omega(U_1),Hor_{y\cdot g}^\omega(U_2)\right)\right]_{\mathcal G}\quad=\quad\left[y\cdot g,-\omega_{y\cdot g}\left(\left[Hor_{y\cdot g}^\omega(U_1),Hor_{y\cdot g}^\omega(U_2)\right]\right)\right]_{\mathcal G}\vspace{0.1cm}\\
\quad\begin{array}{cl}
= & \left[y\cdot g,-\omega_{y\cdot g}\left(\left[(d\Phi)_{(y,g)}\left(Hor_y^\omega(U_1),(d\gamma)_x(u_1)\right),(d\Phi)_{(y,g)}\left(Hor_y^\omega(U_2),(d\gamma)_x(u_2)\right)\right]\right)\right]_{\mathcal G}\vspace{0.1cm}\\
\overset{(\star)}{=} & \left[y\cdot g,-\omega_{y\cdot g}\left((d\Phi)_{(y,g)}\left(\left[Hor_y^\omega(U_1),Hor_y^\omega(U_2)\right],[(d\gamma)_x(u_1),(d\gamma)_x(u_2)]\right)\right)\right]_{\mathcal G}\vspace{0.1cm}\\
= & \left[y\cdot g,-\omega_{y\cdot g}\left((d\Phi)_{(y,g)}\left(\left[Hor_y^\omega(U_1),Hor_y^\omega(U_2)\right],(d\gamma)_x\left([u_1,u_2]\right)\right)\right)\right]_{\mathcal G}\vspace{0.1cm}\\
= & \left[y\cdot g,-Ad_{g^{-1}}\left(\omega_y\left(\left[Hor_y^\omega(U_1),Hor_y^\omega(U_2)\right]\right)\right)\right]_{\mathcal G}\vspace{0.1cm}\\
= & \left[y,-\omega_y\left(\left[Hor_y^\omega(U_1),Hor_y^\omega(U_2)\right]\right)\right]_{\mathcal G}\vspace{0.1cm}\\
= & \left[y,\Omega_y\left(Hor_y^\omega(U_1),Hor_y^\omega(U_2)\right)\right]_{\mathcal G},
\end{array}
\end{array}
\end{equation}
where we have used that $\left[(d\gamma)_x(u_1),(d\gamma)_x(u_2)\right]=(d\gamma)_x\left([u_1,u_2]\right)$.

%%%%%%%%%%%%%%%%%%%%
\section{Examples}

%%%%%%%%%%%%%
\subsection{Standard principal bundles and connections}

Generalized principal connections reduce to usual principal connections on (standard) principal bundles (for exmaple, \cite[Ch.II]{kobayashinomizu1963}, \cite[Ch.III]{michor1993}). Let $\pi_{P,X}$ be a principal $G$-bundle and denote by $R\colon P\times G\to P$ the corresponding right action. We define a fibered action of the trivial Lie group bundle $\mathcal G=X\times G$ on $\pi_{P,X}$ as
\begin{equation*}
\begin{array}{rccl}
\Phi\colon & P\times_X\mathcal G & \longrightarrow & P\\
& (y,(x,\hat g)) & \longmapsto & R_{\hat g}(y)=y\cdot\hat g
\end{array}
\end{equation*}
Note that $P/\mathcal G\simeq X$, so we can regard $\pi_{P,X}$ as a generalized principal bundle with respect to this fibered action. 

Let $\nu_0$ be the trivial connection on $\pi_{\mathcal G,X}$, that is, the one given by
\begin{equation*}
(\nu_0)_{(x,\hat g)}\left(U_x,U_{\hat g}\right)=(0_x,U_{\hat g}),\qquad  (U_x,U_{\hat g})\in T_{(x,\hat g)}\mathcal G=T_x X\oplus T_{\hat g} G
\end{equation*}
In addition, note that the Lie algebra bundle of $\pi_{\mathcal G,X}$ is $\varfrak g=X\times\mathfrak g$, where $\mathfrak g$ is the Lie algebra of $G$. 

\begin{proposition}\label{prop:principalconnection}
Let $\hat\omega\in\Omega^1(P,\mathfrak g)$ and $\omega\in\Omega^1(P,\varfrak g)$ be such that
\begin{equation*}
\omega_y(U_y)=(x,\hat\omega_y(U_y)),\qquad y\in P,\quad U_y\in T_y P,\quad x=\pi_{P,X}(y).
\end{equation*}
Then $\omega$ is a generalized principal connection associated to $\nu_0$ if and only if $\hat\omega$ is a (standard) principal connection.
\end{proposition}

\begin{proof}
Let $(y,\xi)\in P\times_X\varfrak g$ and denote $\xi=(x,\hat\xi)$, $x=\pi_{P,X}(y)$, for some $\hat\xi\in\mathfrak g$. Observe that
\begin{equation*}
\xi^*_y=\left.\frac{d}{dt}\right|_{t=0}y\cdot\exp(t\,\xi)=\left.\frac{d}{dt}\right|_{t=0}y\cdot\left(x,\exp(t\,\hat\xi)\right)=\left.\frac{d}{dt}\right|_{t=0}y\cdot\exp\left(t\,\hat\xi\right)=\hat\xi^\star_y
\end{equation*}
where $\hat\xi^\star_y$ is the infinitesimal generator of $\hat\xi\in\mathfrak g$ at $y$ in the sense of principal bundles. Hence, 
\begin{equation*}
\omega_y(\xi^*_y)=\xi\quad\Longleftrightarrow\quad\left(x,\hat\omega_y\left(\hat\xi_y^\star\right)\right)=\left(x,\hat\xi\right)\quad\Longleftrightarrow\quad\hat\omega_y\left(\hat\xi^\star_y\right)=\hat\xi .
\end{equation*}
Subsequently, we only need to check $Ad$-equivariance to conclude.

\begin{enumerate}
    \item[$(\Rightarrow)$]  Suppose that $\omega$ is a generalized principal connection on $\pi_{P,X}$ associated to $\nu_0$. Let $y\in P$, $U_y\in T_y P$ and $\hat g\in G$, and denote $g=(x,\hat g)$, $x=\pi_{P,X}(y)$. Pick $\alpha\colon(-\epsilon,\epsilon)\rightarrow P$ such that $\alpha(0)=y$ and $\alpha'(0)=U_y$, we have
    \begin{equation*}
    (dR_{\hat g})_y(U_y)=\left.\frac{d}{dt}\right|_{t=0} R_{\hat g}\left(\alpha(t)\right)=\left.\frac{d}{dt}\right|_{t=0}\Phi\left(\alpha(t),\beta(t)\right)=(d\Phi)_{\left(\alpha(0),\beta(0)\right)}\left(\alpha'(0),\beta'(0)\right)
    \end{equation*}
    where $\beta\colon(-\epsilon,\epsilon)\rightarrow\mathcal G$ is given by $\beta(t)=\left((\pi_{P,X}\circ\alpha)(t),\hat g\right)$. Note that $\beta(0)=(x,\hat g)=g$ and $\beta'(0)=(U_x,0_{\hat g})$ for $U_x=(d\pi_{P,X})_y(U_y)\in T_x X$. As a result,
    \begin{equation*}
    (dR_{\hat g})_y(U_y)=(d\Phi)_{\left(y,g\right)}\left(U_y,U_g\right),\qquad U_g=(U_x,0_{\hat g})
    \end{equation*}
    Now, property $(ii)$ of Definition \ref{def:connectionform} ensures that $\hat\omega$ is a principal connection on $\pi_{P,X}$:
    \begin{equation*}
    \begin{array}{ccl}
    (x,\hat\omega_{y\cdot\hat g}\left((dR_{\hat g})_y(U_y)\right) & = & \omega_{y\cdot g}\left((d\Phi)_{\left(y,g\right)}\left(U_y,U_g\right)\right)\\
    & = & Ad_{g^{-1}}\big(\omega_y(U_y)+(\nu_0)_g(U_x,0_{\hat g})\big)\\
    & = & Ad_{g^{-1}}\left(x,\hat\omega_y(U_y)\right)\\
    & = & \left(x,Ad_{\hat g^{-1}}\left(\hat\omega_y(U_y)\right)\right)
    \end{array}
    \end{equation*}
    
    \item[$(\Leftarrow)$]  Suppose that $\hat\omega$ is a principal connection on $\pi_{P,X}$. Let $(y,g)\in P\times_X\mathcal G$ with $g=(x,\hat g)$, and $(U_y,U_g)\in T_y P\times_{T_x X} T_g\mathcal G$ with $U_g=(U_x,U_{\hat g})\in T_x X\oplus T_{\hat g} G$. Observe that
    \begin{equation*}
    \begin{array}{ccl}
    (d\Phi)_{(y,g)}\left(U_y,U_g\right) & = & (dR)_{(y,\hat g)}\left(U_y,U_{\hat g}\right)\\
    & = & (dR_{\hat g})_y(U_y)+(d\phi_y)_{\hat g}(U_{\hat g})\\
    & = & (dR_{\hat g})_y(U_y)+(d\phi_{y\cdot\hat g})_1\left((dL_{\hat g^{-1}})_{\hat g}(U_{\hat g})\right)\\
    & = & (dR_{\hat g})_y(U_y)+(dL_{\hat g^{-1}})_{\hat g}(U_{\hat g})_{y\cdot\hat g}^\star
    \end{array}
    \end{equation*}
    where $\phi_y\colon G\to P$ is given by $\phi_y(\hat g)=y\cdot\hat g$ and $L_{\hat g}\colon G\rightarrow G$ is the left multiplication by $\hat g$. Using this we conclude
    \begin{equation*}
    \begin{array}{ccl}
    \omega_{y\cdot g}\left((d\Phi)_{(y,g)}(U_y,U_g)\right) & = & \left(x,\hat\omega_{y\cdot\hat g}\left((dR_{\hat g})_y(U_y)+(dL_{\hat g^{-1}})_{\hat g}(U_{\hat g})_{y\cdot\hat g}^\star\right)\right)\vspace{0.1cm}\\
    & = & \left(x,Ad_{\hat g^{-1}}\left(\hat\omega_y(U_y)\right)+(dL_{\hat g^{-1}})_{\hat g}(U_{\hat g})\right)\vspace{0.1cm}\\
    & = & Ad_{g^{-1}}\left(\omega_y(U_y)\right)+(dL_{g^{-1}})_g\left((\nu_0)_g(U_g)\right)\vspace{0.1cm}\\
    & = & Ad_{g^{-1}}\left(\omega_y(U_y)+(\nu_0)_g(U_g)\right)
    \end{array}
    \end{equation*}
\end{enumerate}
\end{proof}

%%%%%%%%%%%%%
\subsection{Affine bundles and connections}\label{affine}

Generalized principal connections on an affine bundle are just affine connections (see \cite{kobayashinomizu1963}, \cite{sardanashvily2013}) associated to linear connections %\cite{darling1994}
on the modelling vector bundle. A vector bundle $\pi_{\overline E,X}$ is an abelian Lie group bundle with the additive structure. A Lie group connection $\nu\colon T\overline E\rightarrow\overline E$ is just a linear connection, since it respects this additive structure. If we consider linear bundle coordinates $(x^\mu,v^A)$ corresponding to a basis of local sections $\{e_A\colon 1\leq A\leq m\}$ of $\pi_{\overline E,X}$, then we may write
\begin{equation*}
\nu(x^\mu,v^A)=\nu_{\mu,B}^A(x^\mu)\,v^B\,dx^\mu\otimes e_A+dv^A\otimes e_A,
\end{equation*}
for some (local) functions $\nu_{\mu,A}^B\in C^\infty(X)$, $1\leq\mu\leq n$, $1\leq A,B\leq m$.

Now consider an affine bundle $\pi_{E,X}$ modelled on $\pi_{\overline E,X}$. We have the following fibered action
\begin{equation*}
\begin{array}{ccl}
E\times_X\overline E & \longrightarrow & E\\
(y,v) & \longmapsto & y+v
\end{array}
\end{equation*}
This action is clearly free and proper, and it satisfies that $E/\overline E\simeq X$, so we can regard $\pi_{E,X}$ as a generalized principal bundle. In the following result we see that generalized principal connections for this action are just affine connections.

\begin{proposition}
Let $\omega\in\Omega^1\left(E,\overline E\right)$ and $\nu$ be a linear connection on $\pi_{\overline E,X}$. Then $\omega$ is an affine connection on $\pi_{E,X}$ with $\nu$ as underlying linear connection if and only if $\omega$ is a generalized principal connection on $\pi_{E,X}$ associated to $\nu$.
\end{proposition}

\begin{proof}
Let $(x^\mu,y^A)$ be affine bundle coordinates for $\pi_{E,X}$ associated to $(x^\mu,v^A)$. First, we suppose that $\omega$ is an affine connection with $\nu$ as underlying linear connection. The local expression of $\omega$ is
\begin{equation*}
\omega(x^\mu,y^A)=\left(\nu_{\mu,B}^A(x^\mu)y^B+\Gamma_\mu^A(x^\mu)\right)dx^\mu\otimes e_A+dy^A\otimes e_A
\end{equation*}
for some (local) functions $\Gamma_\mu^A\in C^\infty(X)$, $1\leq\mu\leq n$, $1\leq A\leq m$. It is straightforward to check the the equivariance condition from Definition \ref{def:connectionform}, that for this case reads $\omega_{y+v}=\omega_y+\nu_v$ for each $(y,v)\in E\times_X\overline E$ (observe that the adjoint representation is trivial, since the group is abelian).

Conversely, suppose that $\omega$ is a generalized principal connection associated to $\nu$. Locally, $\omega$ will be given by
\begin{equation*}
\omega(x^\mu,y^A)=\omega_\mu^A(x^\mu,y^A)\,dx^\mu\otimes e_A+dy^A\otimes e_A
\end{equation*}
for some (local) functions $\omega_\mu^A\in C^\infty(E)$, $1\leq\mu\leq n$, $1\leq A\leq m$. It must satisfy the equivariance condition, that is, $\omega(x^\mu,y^A+v^A)=\omega(x^\mu,y^A)+\nu(x^\mu,v^A)$ for every $(x^\mu,y^A,v^A)\in E\times_X\overline E$. This gives
\begin{equation*}
\begin{array}{ccl}
\omega_\mu^A(x^\mu,y^A+v^A)=\omega_\mu^A(x^\mu,y^A)+\nu_{\mu,B}^A(x^\mu)v^B,\qquad 1\leq\mu\leq n,\quad 1\leq A\leq m.
\end{array}
\end{equation*}
Thence, $\omega_\mu^A(x^\mu,y^A)=\Gamma_\mu^A(x^\mu)+\nu_{\mu,B}^A(x^\mu)y^B$ for $1\leq\mu\leq n$ and $1\leq A\leq m$, where $\Gamma_\mu^A(x^\mu)=\omega_\mu^A(x^\mu,0)$. This is the local expression of an affine connection with $\nu$ as underlying linear connection.
\end{proof}

%%%%%%%%%%%%%
\subsection{Gauge transformations}

Let $\pi_{P,X}\colon P\to X$ be a (standard) principal bundle with structure group $G$. Recall that gauge transformations on $\pi_{P,X}$ are in a bijective correspondence with sections of the Lie group bundle 
\begin{equation*}
Ad(P)=(P\times G)/G\longrightarrow X,    
\end{equation*}
where the right action of $G$ on $P\times G$ is given by $(p,g)\cdot h=(p\cdot h,h^{-1} g h)$ for each $p\in P$ and $g,h\in G$. This provides a fibered action of $\pi_{Ad(P),X}$ on $\pi_{P,X}$. The adjoint bundle $Ad(P)$ is a main instance of Lie group bundle and, in particular, the question of when an arbitrary Lie group bundle $\mathcal{G}\to X$ is the adjoint bundle of a (standard) principal bundle has interesting topological consequences (for example, see \cite{Mackenzie1989}).

Principal connections on $\pi_{P,X}$ are in a bijective correspondence with sections of the bundle of connections,
\begin{equation*}
C(P)=(J^1 P)/G\longrightarrow X,
\end{equation*}
which is an affine bundle modelled on $T^*X\otimes ad(P)\to X$, of covectors taking values in the adjoint bundle $ad(P)=(P\times\mathfrak g)/G$, where $\mathfrak g$ is the Lie algebra of $G$ and the action of $G$ on $\mathfrak g$ is given by the adjoint representation. The group of gauge transformations acts on connections. This action is fiberwisely expressed as 
\begin{eqnarray*}
J^1 Ad(P)\times_X C(P)&\longrightarrow& C(P)\\
(j^1_x\gamma,[j^1_x s])&\mapsto & [j^1_x (\gamma\circ s)].
\end{eqnarray*}
The 1-jet lift of this action is
\begin{equation*}
J^1\left(J^1 Ad(P)\right)\times_X J^1 C(P)\longrightarrow J^1 C(P).
\end{equation*}
This is again a (left) fibered action, but it is not free. Nevertheless, we may consider the following Lie group subbundle of $J^1\left(J^1 Ad(P)\right)$,
\begin{equation*}
J_0^2 Ad(P)=\left\{j_x^2\gamma\in J^2 Ad(P)\mid \gamma(x)=1_{Ad(P)_x}\right\}\subset J^2 Ad(P)\subset J^1\left(J^1 Ad(P)\right).
\end{equation*}
The restriction of the action to this Lie group subbundle is free and proper, making
\[
J^1C(P) \longrightarrow J^1C(P)/J^2_0Ad(P)
\]
a generalized principal bundle that is not a (standard) principal bundle. Furthermore, the quotient is isomorphic to the curvature bundle of $\pi_{P,X}$,
\begin{equation*}
\begin{array}{ccc}
J^1 C(P)/J_0^2 Ad(P) & \widetilde{\longrightarrow} & \bigwedge^2 T^*X\otimes ad(P)\vspace{0.2cm}\\
\left[j_x^1 A\right]_{J_0^2 Ad(P)} & \longmapsto & (F_A)_x,
\end{array}
\end{equation*}
where we denote by $F_A\in\Omega^2(X,ad(P))$ the reduced curvature corresponding to a principal connection $A\in\Omega^1(P,\mathfrak g)$. This is (the main part of) the geometric statement of the well-known Utiyama theorem (for example, see \cite[\S5]{garcia1977}).

\subsubsection{An example of generalized principal connection}

Within the framework of gauge field theories analyzed above, we now present a simple example of a generalized principal connection. Even though one initially would like to have such connection on $C(P)=(J^1P)/G$, the action of $J^1 Ad(P)$ on $C(P)$ is not free. To remedy this, we may consider the action on $J^1 P$ instead,
\begin{equation*}
\begin{array}{rccc}
\Phi\colon & J^1 Ad(P)\times_X J^1 P & \longrightarrow & J^1 P\\
& (j^1_x\varphi, j^1_x s) & \longmapsto & j^1_x(\varphi \circ s).
\end{array}
\end{equation*}
that is free and and proper. Actually, this construction is conceptually close to the case of affine connections (see \S \ref{affine} above), since $J^1 Ad(P)$ acts transitively on $J^1P$.

For the sake of brevity, in the following, we work in a trivializing chart of $\pi_{P,X}$, i.e. we pick $\mathcal{U}\subset X$ with $\mathcal{U}\simeq \mathbb{R}^n$, so that we can write $P|_{\mathcal{U}}=\mathcal U\times G$, $Ad(P)|_{\mathcal{U}}=\mathcal U\times G$ and $ad(P)|_{\mathcal{U}}=\mathcal U\times\mathfrak g$. For the sake of simplicity, we write $\mathcal U=X$. Furthermore, consider the right trivialization $TG\simeq G\times\mathfrak g$, which is given by $U_g\mapsto\left(g,(dR_{g^{-1}})_g(U_g)\right)$. Using this, we have
\begin{equation*}
J^1 P=J^1 Ad(P)=G\ltimes(T^*X\otimes\mathfrak g),\qquad J^1 ad(P)=\mathfrak g\oplus(T^*X\otimes\mathfrak g),    
\end{equation*}
where the (fiberwise) semidirect product $\ltimes$ is given by the adjoint representation, that is,
\begin{equation*}
\left(g,\xi_x\right)\left(g',\xi_x'\right)=\left(gg',\xi_x+Ad_g\circ\xi_x'\right),\qquad\left(g,\xi_x\right),\left(g',\xi_x'\right)\in J_x^1 Ad(P)=G\ltimes(T_x^* X\otimes\mathfrak g).
\end{equation*}

\begin{remark}
It is easy to check that 
\begin{equation*}
(g,\xi_x)^{-1}=\left(g^{-1},-Ad_{g^{-1}}\circ\xi_x\right),\qquad(g,\xi_x)\in J^1 Ad(P).
\end{equation*}
In addition, the adjoint representation of $J^1 Ad(P)$ is given by
\begin{equation*}
Ad_{(g,\xi_x)}(\eta,\phi_x)=\left(Ad_g(\eta),Ad_g\circ\phi_x-[Ad_g(\eta),\xi_x]\right)
\end{equation*}
for every $(g,\xi_x)\in J^1 Ad(P)$ and $(\eta,\phi_x)\in J^1 ad(P)$, where $[\cdot,\cdot]$ is the Lie bracket of $\mathfrak g$.
\end{remark}

Under this trivialization, the (left) fibered action is given by
\begin{equation}\label{eq:accionfibradayangmills}
\left(g,\xi_x\right)\cdot\left(h,A_x\right)=\left(gh,Ad_g\circ A_x+\xi_x\right),\qquad \left(\left(g,\xi_x\right),\left(h,A_x\right)\right)\in J^1 Ad(P)\times_X J^1 P.
\end{equation}

\begin{lemma}\label{lemma:secciongauge}
If $P=X\times G$, we have the identification $$J^1(J^1 Ad(P))=G\ltimes\left((T^*X\otimes\mathfrak g)\oplus(T^*X\otimes\mathfrak g)\oplus(T^*X\otimes T^*X\otimes\mathfrak g)\right).$$ Then, the jet section 
$\hat\nu\in\Gamma\left(\pi_{J^1\left(J^1 Ad(P)\right),J^1 Ad(P)}\right)$ defined as
\begin{equation*}
\hat\nu\left(g,\xi_x\right)=\left(g,\xi_x,\xi_x,0_x\right),\qquad\left(g,\xi_x\right)\in J^1 Ad(P)  
\end{equation*}
is a Lie group bundle connection on $\pi_{J^1Ad(P),X}$.
\end{lemma}

\begin{proof}
It is clear that $(\hat\nu\circ 1)(x)=\hat\nu(1,0_x)=(1,0_x,0_x,0_x)=(d1)_x$ for each $x\in X$. Besides, we have
\begin{equation*}
\begin{array}{ccl}
\hat\nu\left(\left(g,\xi_x\right)\left(g',\xi_x'\right)\right) & = & \hat\nu\left(gg',\xi_x+Ad_g\circ\xi_x'\right)\\
& = & \left(gg',\xi_x+Ad_g\circ\xi_x',\xi_x+Ad_g\circ\xi_x',0_x\right)\\
& = & \left(g,\xi_x,\xi_x,0_x\right)\left(g',\xi_x',\xi_x',0_x\right)\\
& = & \hat\nu\left(g,\xi_x\right)\hat\nu\left(g',\xi_x'\right)
\end{array}
\end{equation*}
and we conclude by Proposition \ref{prop:liegroupconnectionjet}.
\end{proof}

Observe that the Lie algebra bundle of $J^1 Ad(P)$ is $J^1 ad(P)$. Theorem \ref{theorem:existencia} ensures that there exists a generalized principal connection $\omega\in\Omega^1(J^1 P,J^1 ad(P))$ associated to $\nu$. 

\begin{proposition}
If $P=X\times G$, we have the identification $$J^1(J^1 P)=G\ltimes\left((T^*X\otimes\mathfrak g)\oplus(T^*X\otimes\mathfrak g)\oplus(T^*X\otimes T^*X\otimes\mathfrak g)\right).$$ An Ehresmann connection $\omega\in\Omega^1(J^1 P,J^1 ad(P))$ with corresponding jet section $\hat\omega\in\Gamma(\pi_{J^1(J^1 P),J^1 P})$ is a generalized principal connection on $\pi_{J^1 P,X}$ associated to $\nu$ if and only if
\begin{equation*}
\hat\omega(h,A_x)=\left(h,A_x,Ad_h\circ f(x)+A_x,Ad_h\circ g(x)\right),\qquad(h,A_x)\in J^1 P,
\end{equation*}
for some sections $f\in\Gamma(\pi_{T^*X\otimes\mathfrak g,X})$ and $g\in\Gamma(\pi_{T^*X\otimes T^*X\otimes\mathfrak g,X})$.
\end{proposition}

\begin{proof}
Observe that we may write
\begin{equation*}
\hat\omega(h,A_x)=\left(h,A_x,\varpi(h,A_x),\tilde\varpi(h,A_x)\right),\qquad (h,A_x)\in J^1 P,
\end{equation*}
for some $\varpi\colon J^1 P\to T^*X\otimes\mathfrak g$ and $\tilde\varpi\colon J^1 P\to T^*X\otimes T^*X\otimes\mathfrak g$. On the other hand, it can be seen that the jet extension of the action \eqref{eq:accionfibradayangmills} is given by
\begin{equation*}
\begin{array}{l}
\left(g,\xi_x,\eta_x,\phi_x\right)\cdot\left(h,A_x,\chi_x,\alpha_x\right)=\\
\left(gh,Ad_g\circ A_x+\xi_x;Ad_g\circ\chi_x+\eta_x,Ad_g\circ\alpha_x+\phi_x+[\eta_x,Ad_g\circ A_x]\right)
\end{array}
\end{equation*}
for each $(g,\xi_x,\eta_x,\phi_x)\in J^1(J^1Ad(P))$ and $(h,A_x,\chi_x,\alpha_x)\in J^1(J^1 P)$. By choosing $(h,A_x)=(1,0_x)$ we get
\begin{equation*}
\hat\omega\left((g,\xi_x)\cdot(1,0_x)\right)=\hat\omega\left(g,\xi_x\right)=\left(g,\xi_x,\varpi(g,\xi_x),\tilde\varpi(g,\xi_x)\right).
\end{equation*}
Likewise,
\begin{equation*}
\begin{array}{ccl}
\hat\nu(g,\xi_x)\cdot\hat\omega(1,0_x) & = & \left(g,\xi_x,\xi_x,0_x\right)\cdot\left(1,0_x,\varpi(1,0_x),\tilde\varpi(1,0_x)\right)\\
& = & \left(g,\xi_x,Ad_g\circ\varpi(1,0_x)+\xi_x,Ad_g\circ\tilde\varpi(1,0_x)\right).
\end{array}
\end{equation*}
Proposition \ref{prop:jetsectionomega} ensures that $\omega$ is a generalized principal connection on $\pi_{J^1 P,X}$ associated to $\nu$ if and only if
\begin{equation*}
\varpi(g,\xi_x)=Ad_g\circ\varpi(1,0_x)+\xi_x,\qquad\tilde\varpi(g,\xi_x)=Ad_g\circ\tilde\varpi(1,0_x)
\end{equation*}
for each $(g,\xi_x)\in J^1 P$. By denoting $f(x)=\varpi(1,0_x)$ and $g(x)=\tilde\varpi(1,0_x)$ for each $x\in X$ we conclude.
\end{proof}

\paragraph{Acknowledgments.}

MCL has been partially supported by grant no. PGC2018-098321-B-I00, Ministerio de Ciencia e Innovación, Spain; and grant no. SA090G19, Consejería de Educación, Junta de Castilla y León, Spain.

ARA has been supported by a FPU grant from the Spanish Ministry of Science, Innovation and Universities (MICIU).

%%%%%%%%%%%%%%%%%%%%%%%%%%%%%%%%%%%%%%%%%%%%%%%%%%
\bibliographystyle{plain}
\bibliography{biblio_reduction.bib}

\end{document}